\documentclass[11pt]{article}

\usepackage[margin=1in]{geometry}

\usepackage{amsmath,amscd}
\usepackage{amssymb}
\usepackage{amsthm}
\usepackage{comment}
\usepackage{graphicx}
\usepackage{subcaption} 
\usepackage{epstopdf} 
\usepackage{enumitem}
\usepackage{mathrsfs,mathtools}
\usepackage{cite}
\usepackage{authblk}
\usepackage{bm}
\usepackage{hyperref}
\usepackage{xcolor}

\newtheorem{thm}{Theorem}[section]

\newtheorem*{proposition*}{Proposition}
\numberwithin{equation}{section}

\allowdisplaybreaks

\title{Parameter Identification in Reaction-Diffusion-Chemotaxis Systems from Single Turing Pattern Amplitudes}

\author[1,*]{Yuhan Li}
\author[1,$\dagger$]{Hongyu Liu}
\author[2,$\natural$]{Catharine W. K. Lo}
\affil[1]{Department of Mathematics, City University of Hong Kong, Hong Kong, P.R. China}
\affil[2]{School of Mathematical Sciences, Shenzhen University, Shenzhen, 518061, P.R. China}
\affil[*]{yuhli2-c@my.cityu.edu.hk}
\affil[$\dagger$]{hongyliu@cityu.edu.hk, hongyu.liuip@gmail.com}
\affil[$\natural$]{cwklo@szu.edu.cn, catharinelowk@gmail.com}
\date{}

\begin{document}
\maketitle
\makeatletter
\newcommand{\rmnum}[1]{\romannumeral #1}
\newcommand{\Rmnum}[1]{\expandafter\@slowromancap\romannumeral #1@}
\makeatother

\begin{abstract}
Turing patterns play a fundamental role in morphogenesis and population dynamics, encoding key information about the underlying biological mechanisms. Traditional inverse problems, however, have largely relied on non-biological data such as boundary measurements, neglecting the rich information embedded in the patterns themselves. In this work, we introduce a new research direction that directly leverages physical observables from nature -- the amplitudes of Turing patterns -- to achieve parameter identification. Specifically, motivated by stripe-like patterns commonly observed in biological systems, we consider two one-dimensional models with prescribed nonlinearities: one with density-dependent chemotaxis ($\chi(n,c)=\chi_0 n$) and one with ratio-dependent chemotaxis ($\chi(n,c)=\chi_0 n/c$). We present a mathematical framework that uses the spatial amplitude profile of a stationary pattern to recover the system parameters, including the wavelength, diffusion constants, and chemotactic and kinetic coefficients. This amplitude-based approach establishes a biologically grounded and mathematically rigorous paradigm for reverse-engineering pattern formation mechanisms across diverse biological systems.

The core mathematical contribution is a proof-of-concept demonstration that, under a Fourier truncation at mode $M=3$ (justified by the rapid decay of higher harmonics near the Turing instability threshold), the amplitude data of a stationary pattern can be used to reconstruct the parameter combinations $\{d_n k^2, d_c k^2, \chi_0 k^2, r\}$ for both models. These combinations uniquely determine the ratios $d_c/d_n$ and $\chi_0/d_n$, while the individual parameters are recovered up to an overall scale. We emphasize that this is a purely theoretical work; numerical validation and stability analysis are not included but constitute important directions for future research. The framework is extensible to other models with suitable modifications.
~\\\\
\textbf{Keywords:} Inverse reaction-diffusion equations; Turing patterns; Turing instability; periodic solutions; sinusoidal form.~\\
\textbf{2020 Mathematics Subject Classification:} 35R30, 35B10, 35B36, 35K10, 35K55, 35K57, 35Q92, 92-10, 92C15, 92C37, 92C70, 92D25
\end{abstract}

\section{Introduction}  \label{introduction}

\subsection{Background and motivation}\label{sec:Turing_background}
Nature displays an astonishing variety of organized spatial structures~\cite{Cantrell2020,gupta2009linear, kaushal2011human}, from animal coat patterns to bacterial colonies. These recurring morphologies are not accidental but arise from reproducible developmental or ecological processes, suggesting universal underlying mechanisms. A powerful explanation was proposed by Alan Turing in 1952 \cite{Turing90}: the interaction of two morphogens -- an activator and an inhibitor -- through reaction and diffusion can spontaneously break symmetry and generate stable heterogeneous patterns. This mechanism, now known as Turing instability, occurs when the inhibitor diffuses faster than the activator, amplifying small fluctuations into periodic concentration profiles that manifest as stripes, spots, or labyrinths.

\begin{figure}[htp]
\begin{center}
\begin{subfigure}{0.45\textwidth}
        \centering
        \includegraphics[width=\linewidth]{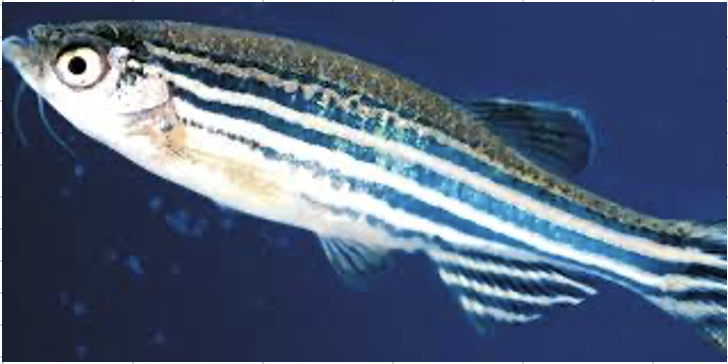}
        \caption{Wild-type adult zebrafish stripe pattern. \textbf{Source:} \cite{leiden_university_2015}} 
        \label{fig:5.1a}
    \end{subfigure}
    \quad
    \begin{subfigure}{0.45\textwidth}
        \centering
        \includegraphics[width=\linewidth]{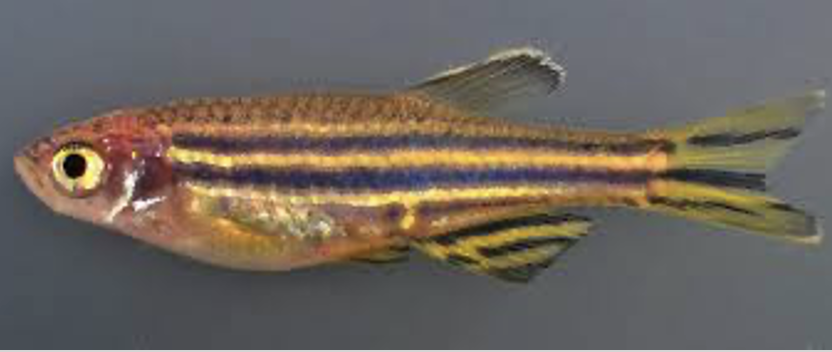}
        \caption{Variant zebrafish stripe pattern. \textbf{Source:} \cite{eurekalert_news}} 
        \label{fig:5.1b}
    \end{subfigure}
  \caption{Adult pigmentation patterns of zebrafish. The pattern is formed by coordinated arrangement of melanophores, xanthophores, and iridophores, representing a typical vertebrate biological pattern formation system.}
    \label{fig:5.1} 
  \end{center}
\end{figure}

\begin{figure}[htp]
    \begin{center} 
    \begin{subfigure}{0.45\textwidth}
        \centering
        \includegraphics[width=\linewidth]{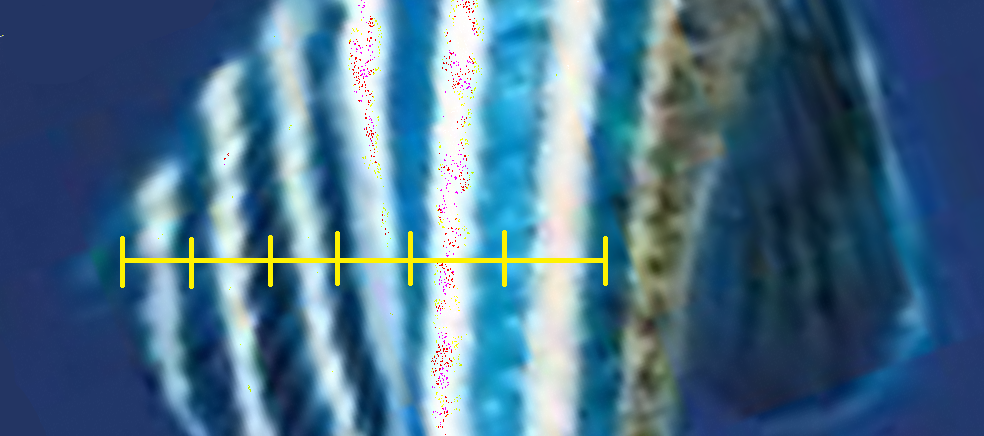}
        \caption{One dimensional periodic profile of Figure 1a. } 
        \label{fig:5.3a}
    \end{subfigure}
    \quad
    \begin{subfigure}{0.45\textwidth}
        \centering
        \includegraphics[width=\linewidth]{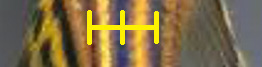}
        \caption{One dimensional periodic profile of Figure 1b. } 
        \label{fig:5.3b}
    \end{subfigure}
    \caption{1D transects of adult pigmentation patterns of zebrafish.}
    \label{fig:5.2} 
    \end{center}
\end{figure}

Mathematically, such dynamics are often modeled by reaction-diffusion-advection systems, which incorporate directed motion such as chemotaxis \cite{painter2011spatio,ma2020stationary}. These models have successfully described pattern formation in skin pigmentation, bacterial aggregation, vascular networks, and tumor invasion \cite{Bard1981model, Murray1, Murray2, Steeves1989patterns, Ball2015forging}. 
Figure \ref{fig:5.1} illustrates a classic biological example of stripe patterns in zebrafish, demonstrating the kind of periodic spatial structures that motivate our inverse problem framework.
While much research has focused on predicting patterns from given parameters (the forward problem), the inverse problem of reconstructing the underlying mechanism from observed patterns has received less attention, especially in biologically realistic settings. In practice, researchers often have access to spatial or spatio-temporal pattern data, such as concentration profiles, color intensities, or density maps. Among these, the spatial amplitude of a pattern, which measures deviation from uniformity, is often experimentally accessible and informative: from a one-dimensional transect, one can compute Fourier coefficients $\alpha_i = \frac{2}{L}\int_0^L n(x)\cos(ikx)dx$ (see Section \ref{M_analysis} for details). 
Such transects are biologically relevant for patterns with straight stripe geometries -- common in zebrafish pigmentation, bacterial colony formations, and certain plant leaf venation patterns -- where the pattern is approximately periodic in one direction and translationally invariant in the transverse direction (see Figures \ref{fig:5.2}(a) and \ref{fig:5.2}(b) for schematic illustrations).
This raises a key question: \emph{Can the amplitude information embedded in a stationary Turing pattern be used to determine the system parameters for specific reaction-diffusion-chemotaxis models?}

In this work, we provide an affirmative answer for a class of one-dimensional reaction-diffusion-chemotaxis models with density-dependent feedback. We show theoretically that the spatial amplitude profile of a single pattern contains sufficient information to uniquely determine certain parameter combinations, from which the individual parameters can be recovered up to an overall scale. Our approach establishes a direct link between observable pattern form and generative mechanism for two specific models with prescribed nonlinearities, and provides a rigorous inverse-problem framework for these cases. A detailed introduction to the methodology and its novelty is given in Section \ref{sec:Turing_Novelty}.

\subsection{Model formulation}\label{sec:Turing_ModelF}
The overarching framework for our study is a class of reaction-diffusion-advection systems designed to model the population dynamics of motile micro-organisms (e.g., bacteria or amoeba) in response to a chemotactic agent \cite{Bucur2024,Hillen2009user,Keller1970initiation}. The general system is given by:
\begin{equation}\label{originial_system}
     \begin{dcases}
         \partial_t n=d_n\Delta n-\nabla\cdot(\chi(n,c)\nabla c)+f(n,c),\\
         \partial_t c=d_c\Delta c+g(n,c),
     \end{dcases}
 \end{equation}
where the variable $n$ denotes the micro-organism density, and $c$ represents the chemotactic agent concentration. The parameter $d_n$ characterizes the diffusion of the micro-organisms as a result of random cellular movement, $d_c$ signifies the chemotactic agent's diffusion constant, 

This model serves as a suitable framework for investigating the mechanisms underlying pattern formation in populations of motile micro-organisms. By considering various forms for the functions $\chi(n,c)$, $f(n,c)$, and $g(n,c)$, we can adapt \eqref{originial_system} into models with diverse physical applications.

A foundational special case is the following cross-diffusion model \cite{Keller1970initiation}, obtained by setting $\chi(n,c)\equiv \chi_0$ and $f(n,c)\equiv0$:
\begin{equation}\label{KS_model}
    \begin{dcases}
        \partial_t n=d_n \Delta n-\chi_0\Delta c,\\
        \partial_t c=d_c \Delta c+hn-pc.
    \end{dcases}
\end{equation}
Originally proposed to elucidate the movement and clustering behavior of \textit{Dictyostelium discoideum} amoeba, this model demonstrates how chemotaxis can drive pattern formation. The homogeneous steady state $(n_0,c_0) = (n_0, n_0 h/p)$ becomes unstable when the condition
\[\frac{\chi_0 p}{d_n h}+\frac{n_0 h}{p}>1\]
is met, leading to cell aggregation and the emergence of stationary periodic patterns through Turing instability.

However, constant chemotactic sensitivities often fail to capture the biological reality where cellular response typically depends on local conditions. In many micro-organisms, chemotactic sensitivity intensifies with cell density, enhancing collective movement toward attractants and creating positive feedback loops that promote aggregation and clustering. This behavior is evident in various natural systems such as bacterial swarming and immune cell responses.

To better represent these biological complexities, we focus on models with density-dependent chemotactic sensitivities. Specifically, we examine two physically relevant cases: $\chi(n,c) = \chi_0 n$ and $\chi(n,c) = \chi_0 n/c$, both of which can produce diverse spatial structures including stripes and spots.

The biological relevance of such approaches is exemplified by \emph{Myxococcus xanthus}, a soil-dwelling bacterium known for its complex social behaviors and ability to form structured fruiting bodies under nutrient scarcity \cite{Black2004myxococcus, Berleman2008predataxis}. This bacterium exhibits notable chemotactic response to A-signal molecules, and its patterning behavior can be captured by models similar to:
\begin{equation}\label{Mxanthus}
    \begin{dcases}
        \partial_t n=d_n \Delta n+\nabla\cdot(\chi(n)\nabla c),\\
        \partial_t c=d_c \Delta c+hn-pc.
    \end{dcases}
\end{equation}
Here, $h$ denotes the production rate, $p$ signifies the degradation rate, and $\chi(n)$ represents the chemotactic sensitivity function that depends on the bacterial cell density. The pattern formation process typically begins when nutrient scarcity prompts cells to release signaling molecules, initiating positive feedback loops that amplify small density fluctuations and eventually lead to stationary periodic Turing patterns.

In this study, we analyze Turing pattern formation in a one-dimensional domain with no-flux boundary conditions. We examine two specific models derived from \eqref{originial_system} that exhibit biologically realistic instability mechanisms. These nonlinear models with density-dependent chemotaxis demonstrate the robustness of our recovery method across different nonlinearities.

\textbf{Model 1 (Density-Dependent Chemotaxis):}
We first consider the system with $\chi(n,c) = \chi_0 n$, $f(n,c) = rn(1-n)$, and $g(n,c) = n-c$:
\begin{equation}\label{M3_original}
    \begin{cases}
        \partial_t n=d_n \Delta n-\chi_0 \nabla\cdot (n \nabla c)+rn(1-n),\\
        \partial_t c=d_c \Delta c +n-c.
    \end{cases}
\end{equation}
This model describes bacterial populations responding to a self-produced chemical attractant. The logistic growth term $rn(1-n)$ accounts for resource-limited population growth, while the term $n-c$ represents chemical production and degradation. The chemotactic term $-\chi_0 \nabla \cdot (n \nabla c)$ indicates movement toward higher chemical concentrations, creating a density-dependent feedback mechanism.

In the absence of diffusion and chemotaxis, the system exhibits an unstable state at $(0,0)$ and a stable state at $(1,1)$. Introducing diffusion produces traveling wavefronts between these states, while adding sufficient chemotactic influence leads to stationary periodic patterns \cite{Bucur2024}.

\textbf{Model 2 (Ratio-Dependent Chemotaxis):}
We also examine the case where $\chi(n,c) = \chi_0 n/c$, with $f(n,c) = rn(1-n)$ and $g(n,c) = n-c$:
\begin{equation}\label{M6_original}
    \begin{cases}
        \partial_t n=d_n \Delta n-\chi_0 \nabla \cdot(\frac{n}{c} \nabla c)+rn(1-n),\\
        \partial_t c=d_c \Delta c +n-c.
    \end{cases}
\end{equation}
This formulation introduces an additional regulatory mechanism where the chemotactic response depends on the ratio of cell density to chemical concentration. This represents scenarios where cells become less responsive in chemical-saturated environments, leading to different instability conditions and pattern characteristics compared to Model 1. Under specific conditions, this system also generates Turing patterns \cite{Bucur2024}.

Both models serve as prototypical frameworks for our inverse problem: from the Fourier amplitude data of a single stationary periodic pattern, we will uniquely recover the full parameter set $\{d_n, d_c, \chi_0, r, k\}$. The detailed recovery procedures are presented in Sections \ref{Proof_M3} and \ref{Proof_M6}, respectively.

\subsection{From biological observations to mathematical quantities}\label{sec:bio_to_math}

Before presenting our inverse results, we clarify how the mathematical quantities in our framework relate to biologically measurable data. This discussion is meant to establish the practical relevance of our theorems, not to report experimental results.

\subsubsection{What can be measured in a biological Turing pattern}

In a biological system exhibiting Turing patterns (e.g., pigmentation stripes, bacterial colonies, or cell density waves), the following measurements are experimentally accessible:

\begin{itemize}
    \item \textbf{Spatial intensity profiles:} Fluorescence microscopy or brightfield imaging provides a two-dimensional map of cell density or pigment concentration. Along a transect perpendicular to the pattern orientation, one obtains a one-dimensional signal $I(x)$. This is particularly natural for stripe-like patterns, where a transect perpendicular to the stripes captures the full spatial variation, while the transverse direction is translationally invariant.

    \item \textbf{Pattern wavelength $\lambda$:} The distance between successive stripes or spots can be measured directly from the image, yielding the fundamental wavenumber $k = 2\pi/\lambda$.
    
    \item \textbf{Baseline density:} The average intensity $\frac{1}{L}\int_0^L I(x)dx$ corresponds to $\alpha_0$ in our Fourier expansion.
    
    \item \textbf{Harmonic amplitudes:} Using the formulas in Section \ref{Proof_M3}, the coefficients $\alpha_i,\beta_i$ for $i=1,2,3$ can be computed from the intensity profile.
\end{itemize}

\subsubsection{What cannot be directly measured (but can be inferred)}

The following quantities are typically \emph{not} directly observable but are targets of our inverse method:

\begin{itemize}
    \item \textbf{Chemical concentrations $c(x)$:} Signaling molecules (e.g., chemoattractants) often require specialized staining or biosensors. Our framework shows they can be inferred from cell-density pattern amplitudes alone.
    
    \item \textbf{Diffusion coefficients $d_n$, $d_c$:} These quantify random motion of cells and molecules. Measuring them directly requires tracking experiments; our method estimates them from pattern shape.
    
    \item \textbf{Chemotactic sensitivity $\chi_0$:} This parameter describes how strongly cells move up chemical gradients -- difficult to measure directly but encoded in pattern sharpness.
\end{itemize}

\subsubsection{Illustrative measurement pipeline (conceptual)}

While the present work is mathematical, we outline a hypothetical pipeline to clarify how our theorems would be applied:

\begin{enumerate}
    \item \textbf{Image acquisition:} Obtain a greyscale image of the pattern (e.g., fluorescently labeled cells).
    \item \textbf{Transect extraction:} Select a line perpendicular to the stripes; record intensity $I(x)$ at positions $x_1, x_2, \dots, x_N$.
    \item \textbf{Preprocessing:} Subtract background (interstripe baseline) and optionally smooth to reduce noise.
    \item \textbf{Wavelength measurement:} Compute average stripe spacing $\lambda$ from peak positions; set $k = 2\pi/\lambda$.
    \item \textbf{Fourier analysis:} Compute coefficients $\alpha_i = \frac{2}{L}\int_0^L I(x)\cos(ikx)dx$ for $i=0,1,2,3$.
    \item \textbf{Parameter recovery:} Solve the algebraic system (e.g., \eqref{M3_equating3}) to obtain $d_n, d_c, \chi_0, r$.
\end{enumerate}

Each step uses only measurements that are, in principle, obtainable from a single biological image. The mathematical guarantee of Theorems \ref{thm:mainthm_M3} and \ref{thm:mainthm_M6} ensures that the recovered parameters are unique.

\subsection{Pattern formation and inverse problem}\label{sec:Turing_PFIP}

The analysis of Turing patterns in system \eqref{originial_system} encompasses two complementary aspects: first, investigating the mechanisms and conditions enabling periodic pattern formation \cite{Turing90, Keller1970initiation}, and second, characterizing the intrinsic properties of the patterns themselves, such as their wavelength and amplitude \cite{Bucur2024,Turing90}.

The wavelength of a Turing pattern represents the spatial periodicity -- the distance between successive concentration peaks or troughs -- that determines whether the pattern manifests as stripes, spots, or other configurations. This fundamental property is governed by parameters within the reaction-diffusion equations, particularly the diffusion rates of the activator and inhibitor components. The amplitude, conversely, quantifies the variation in morphogen concentrations from their homogeneous steady state, reflecting the intensity of patterning, such as the contrast between stripes or the prominence of spots.

In the forward problem, one typically studies Turing patterns arising in two-variable reaction-diffusion systems with polynomial kinetics \cite{Chen2019non}. Under suitable conditions, the solution near the onset of instability can be approximated as a superposition of cosine profiles. In one dimension, this takes the form:
\begin{equation}\label{1D_solution_form}
    \begin{dcases}
        n(x,t)=n_0+\alpha\cos(kx),\\
        c(x,t)=c_0+\beta\cos(kx),
    \end{dcases}
\end{equation}
where $(n_0,c_0)$ is the homogeneous steady state of the system, $k$ represents the wave number of the most unstable mode identified through linear stability analysis, and $\alpha, \beta$ represent the amplitudes of the Turing pattern in the cell density and chemical concentration, respectively. This representation allows estimation of both wavelength $\lambda = 2\pi/k$ and amplitudes from model simulations.

\begin{figure}[htbp]
    \centering
    \begin{subfigure}{0.45\textwidth}
        \centering
        \includegraphics[width=\linewidth]{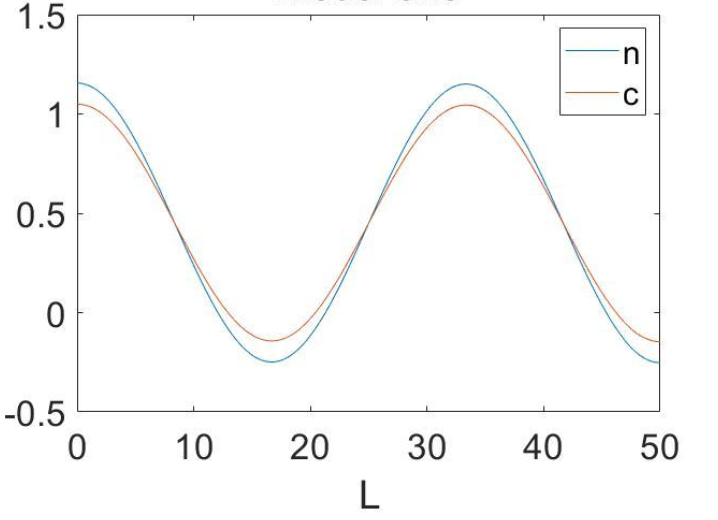}
        \caption{ }  
        \label{pic1}  
    \end{subfigure}
    \hfill
    \begin{subfigure}{0.45\textwidth}
        \centering
        \includegraphics[width=\linewidth]{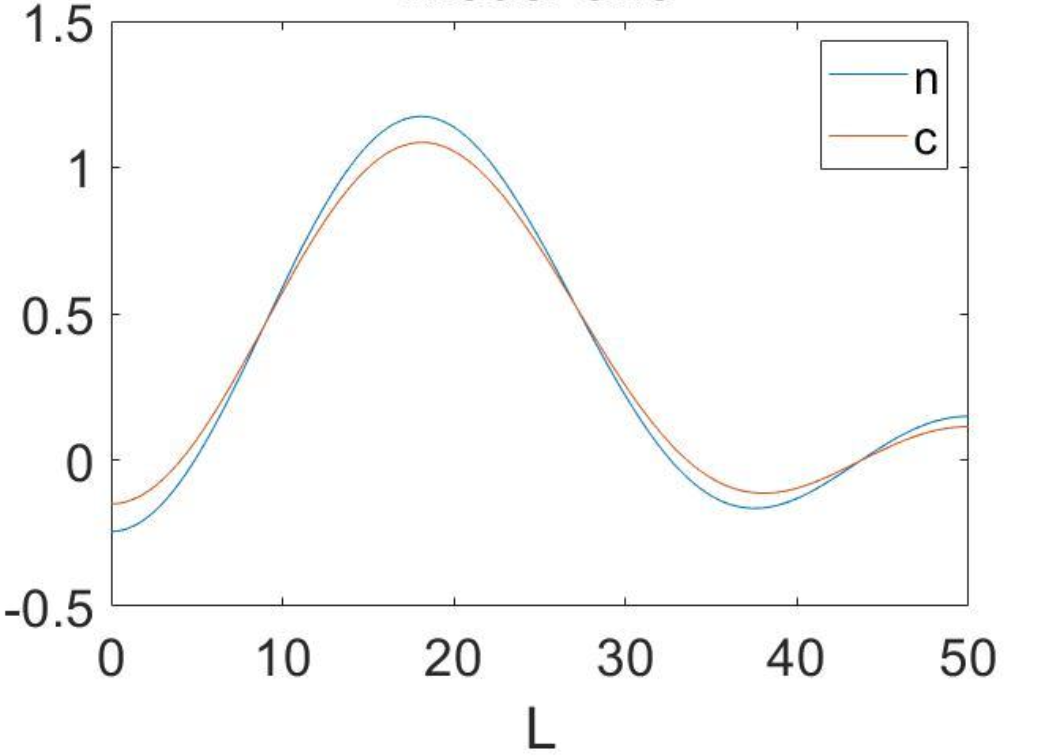}
        \caption{ }  
        \label{pic2}  
    \end{subfigure}
    \caption{}
\end{figure}

The sinusoidal approximation \eqref{1D_solution_form} captures the dominant spatial mode near the instability threshold. In practice, however, Turing patterns may exhibit more structured periodic profiles, as illustrated schematically in Figures \eqref{pic1} and \eqref{pic2}. These examples represent typical spatial distributions of morphogen concentration that arise from reaction-diffusion dynamics, showcasing the kind of periodic amplitude variations that are experimentally observable in biological systems.

Notably, while predicting patterns from model parameters (the forward problem) has been extensively studied, the inverse problem -- recovering parameters from observed patterns -- has received considerably less attention. Yet, in biological contexts, amplitude measurements are often experimentally accessible, for instance, by quantifying color intensity variations in animal coats or contrast in microbial colonies. This observability motivates a fundamental inverse question:

Given complete knowledge of the amplitudes $\alpha$ and $\beta$ in \eqref{1D_solution_form}, is it possible to uniquely determine all unknown parameters in \eqref{originial_system} -- including diffusion coefficients $d_n$, $d_c$, the functional forms of $\chi(n,c)$, $f(n,c)$, $g(n,c)$, and the wavenumber $k$?

To formalize this, we define the inverse problem through the measurement map:
\begin{equation}\label{MeaureMap1}
\mathcal{M}: { \alpha, \beta } \to { d_n, d_c, \chi, f, g, k },
\end{equation}
which encodes the recovery of model parameters from amplitude data. 
Our main theoretical contribution is then summarized as follows:
Given Fourier amplitude measurements 
\begin{equation*}
\mathcal{M}=\{\alpha_i,\beta_i\}_{i=0}^M
\end{equation*}
extracted from a single stationary Turing pattern generated by either Model 1 \eqref{M3_original} or Model 2 \eqref{M6_original}, determine the unknown parameters $\{d_n, d_c, \chi_0, r, k\}$.

Our main results, presented in Sections 3 and 4, show that for both models, under a Fourier truncation at $M=3$ and suitable non-degeneracy conditions on the amplitude data, the measurements uniquely determine the parameter combinations
\begin{equation*}
\{d_n k^2,\; d_c k^2,\; \chi_0 k^2,\; r\}.
\end{equation*}
From these combinations, the ratios $d_c/d_n$ and $\chi_0/d_n$ are uniquely determined, while the individual parameters $\{d_n, d_c, \chi_0, r, k\}$ are recovered up to an overall scale: if $(d_n, d_c, k, \chi_0, r)$ is a solution, then $(\lambda^2 d_n, \lambda^2 d_c, k/\lambda, \lambda^2 \chi_0, r)$ produces the same combinations $d_n k^2$, $d_c k^2$, and $\chi_0 k^2$ for any $\lambda>0$. Additional information, such as knowledge of the domain size or an independent measurement of one parameter, would be required to fix this scale and recover the individual parameters uniquely.

The precise algebraic systems and identifiability statements for each model are given in Theorems \ref{thm:mainthm_M3} and \ref{thm:mainthm_M6} below. We emphasize that the analysis is restricted to one spatial dimension and to the two specific model nonlinearities considered; extensions to higher dimensions or other nonlinearities are not addressed in this work. We emphasize that this is a purely theoretical work; numerical validation and experimental verification are not included but constitute important directions for future research.

\subsection{Discussion of novelty and literature review}\label{sec:Turing_Novelty}
The inverse problem approach developed in this work represents a paradigm shift in both mathematical methodology and biological application. Unlike conventional inverse problems that typically rely on boundary measurements or external data, our framework leverages the intrinsic biological information encoded within Turing patterns themselves to achieve complete parameter identification.

Traditional inverse problems for partial differential equations often employ Dirichlet-to-Neumann (DtN) maps, passive measurements, or single measurement strategies that exploit the mathematical structure through techniques such as high-order variation and complex geometrical optics (CGO) solutions. These approaches, while powerful, often require specialized measurement setups and operate through indirect mathematical reasoning, as seen in studies spanning non-negativity constraints \cite{LLboundary}, Lotka-Volterra systems \cite{LLL1}, chemotaxis models \cite{LLL2,LLL4}, multi-species models \cite{LLL3}, fractional models \cite{ding2024inverse}, and subdomain identification \cite{ding2023determining, ding2024determininginternaltopologicalstructures, LL2025determining, lill2025determining}.

To the best of our knowledge, a sole precedent for parameter recovery from patterns was set by Kazarnikov et al. \cite{Kazarnikov2020statistical}, who used statistical ensemble analysis to recover a single parameter from multiple patterns generated under varying initial conditions. Their method, though innovative, is fundamentally distinct from ours in scope and mechanism: it requires extensive numerical data and statistical correlation, whereas our approach deterministically identifies all model coefficients simultaneously from a single pattern.

Our method represents a radical departure. Rather than relying on boundary measurements, we use the natural observable of a biological system -- the spatially periodic amplitude of a Turing pattern -- as the direct data source. This approach is philosophically inverted: we start from the biological endpoint (the mature pattern) and deduce the generative mechanisms, rather than constraining the problem to mathematically convenient but biologically artificial measurement paradigms.

The mathematical core of our approach rests on a key insight: the exact solution form of a Turing pattern (equation \eqref{1D_solution_form}) contains sufficient information to uniquely determine all underlying parameters. This includes not only the diffusion coefficients $d_n$, $d_c$ and the wavelength $k$, but also the complete functional forms of the chemotactic sensitivity $\chi(n,c)$ and the kinetic terms $f(n,c)$ and $g(n,c)$. This comprehensive identifiability from amplitude data alone bridges a fundamental gap between biological observation and mechanistic modeling.

The two models we analyze serve as proof of concept, but our mathematical framework is general. With suitable computational expansion, the approach can extend to any reaction-diffusion-advection system exhibiting Turing patterns. The recovery mechanism depends only on the pattern-forming capacity of the system, not on specific nonlinearities, representing a new computational paradigm for biological system identification.

Moreover, the general framework \eqref{originial_system} is highly versatile: it reduces to classical predator-prey models ($\chi=0$), encompasses various chemotaxis models, and provides a unified framework for studying pattern formation across biological contexts. Combined with our inverse approach, it opens new avenues for understanding how patterns emerge from underlying mechanisms.

To the best of our knowledge, this constitutes the first comprehensive framework that enables complete parameter identification for Turing pattern-forming systems directly from biological observables. Our work establishes a new research direction that bridges mathematical innovation with biological application, offering a more biologically grounded framework for understanding the complex interplay between genetics, physics, and morphology in developmental processes.

\subsection{Scope and limitations of the present work}\label{sec:limitations}

We emphasize that this is a purely theoretical paper. Our contributions are analytical: we derive algebraic systems that relate Fourier amplitudes to model parameters and argue that these systems are linearly independent. We do not perform numerical simulations, nor do we provide stability estimates for the parameter reconstruction with respect to noise in the measured amplitudes. These constitute important directions for future work.

Furthermore, the analysis is restricted to one spatial dimension and to two specific model nonlinearities. Some biological Turing patterns occur in two or three dimensions, and extending the framework to higher dimensions is meaningful but would require substantial additional work. We note, however, that in many biological systems exhibiting straight stripes -- such as zebrafish pigmentation, certain bacterial colony formations, and some plant patterns -- a one-dimensional transect perpendicular to the stripe orientation captures the essential spatial variation. The present analysis therefore applies to this meaningful class of patterns. Similarly, the generalization to other nonlinearities is not established by our analysis and would need to be verified case by case.

The following sections detail our analysis: Section \ref{sec:Turing_NAresults} presents the nonlinear analysis and main results, with the proof for Model 1 in Section \ref{Proof_M3} and the proof for Model 2 in Section \ref{Proof_M6}. We conclude with remarks and future directions in Section \ref{Crmk}.

\section{Nonlinear analysis and main results} \label{sec:Turing_NAresults}

\subsection{Nonlinear analysis of Turing patterns}  \label{M_analysis} 
As established previously, models such as \eqref{M3_original} and \eqref{M6_original} exhibit stationary periodic Turing patterns, manifesting as highly regular spatial structures including stripes, spots, or labyrinthine forms. A fundamental tool for analysing such patterns is Fourier analysis, originally emphasised by Turing himself in his seminal work \cite{Turing90}. Fourier series provide a natural basis for representing periodic patterns, decomposing complex morphologies into simpler sinusoidal components that illuminate the underlying spatial frequencies and phase relationships governing pattern formation.

To properly characterize the wave profiles in Turing patterns with Fourier series, we express the concentration profiles of the activator and inhibitor substances as a combination of sine and cosine functions. Each term in the series corresponds to a specific frequency, representing different spatial scales of the pattern. By adjusting the coefficients of these sine and cosine terms, we can capture the amplitude and phase of the patterns, allowing us to reconstruct these complex periodic structures.

In this paper, we consider the domain of size $L$ with no-flux boundary conditions, i.e. $\partial_\nu n = \partial_\nu c =0$ at $x=0$ and $x=L$. In this case, the solutions are only represented by the cosine functions, and take the form:
\begin{equation}\label{Fourier_form}
    \begin{cases}
        n(x)=\sum\limits_{i=0}^{M} \alpha_i \cos \left(\dfrac{i\pi x}{L}\right),\\
        c(x)=\sum\limits_{i=0}^{M} \beta_i \cos \left(\dfrac{i\pi x}{L}\right),
    \end{cases}
\end{equation}
where $\alpha_i$ and $\beta_i$ represent the amplitudes of oscillation mode $i$ for the variables $n$ and $c$ respectively. The associated wavelength is given by $k=\pi/L$.  For smooth Turing pattern profiles, the amplitudes $\alpha_i$ and $\beta_i$ quickly tend to zero as $i$ increases. This means that we can truncate the series by considering only the first $M$ terms, with a careful choice of $M$ such that the error is negligible. For further details, refer to the discussion in Section \ref{Proof_M3}.

Given a Turing profile, the amplitudes $\alpha_i$ of $n(x)$ can be computed as follows:
\begin{equation}\label{n_coefficient}
    \alpha_0=\frac{1}{L}\int_0^L n(x)dx, \quad \alpha_i=\frac{2}{L} \int_0^L n(x)\cos(ikx) dx,\ \text{for } i=1,2,\dots.
\end{equation}
Similarly, the amplitudes $\beta_i$ of $c(x)$ can be computed using the following formulae:
\begin{equation}\label{c_coefficient}
    \beta_0=\frac{1}{L}\int_0^L c(x)dx, \quad \beta_i=\frac{2}{L} \int_0^L c(x)\cos(ikx) dx,\ \text{for } i=1,2,\dots.
\end{equation}

The amplitudes $\alpha_i$ and $\beta_i$ quantify the contribution of each cosine component to the overall function. Specifically, $\alpha_0$ and $\beta_0$ represent the average concentrations of the activator and inhibitor over the domain $[0,L]$. This gives a baseline level around which the periodic patterns fluctuate. Other coefficients $\alpha_i$ and $\beta_i$ with $i\geq 1$ correspond to the amplitudes of the cosine functions at different harmonic frequencies $ik$. Each of them quantifies how much of the respective cosine wave contributes to the overall shape of the concentration profiles.

In summary, the amplitudes reveal the predominant spatial frequencies present in the pattern, contributing to a more profound understanding of the Turing structures at play. Through a detailed scrutiny of these coefficients, one can unearth the underlying model and the necessary conditions pivotal for their emergence.

\subsection{Well-posedness of the forward problem}

Before addressing the inverse problem, we briefly discuss the well-posedness of the forward problem, namely, the initial-boundary value problem associated with systems \eqref{M3_original} and \eqref{M6_original} on a bounded domain subject to no-flux boundary conditions. The purpose of this discussion is to specify the class of forward solutions on which the subsequent inverse analysis is based. This ensures that for any given set of parameters \((d_n,d_c,\chi_0,r)\) and non-negative initial data, there exists a unique solution, and that solution depends continuously on the parameters, justifying the use of Fourier analysis on numerically computed patterns.

For Model 1 \eqref{M3_original}, the chemotactic flux is of the standard form
\[
-\chi_0\nabla\cdot(n\nabla c),
\]
while the logistic term
\[
rn(1-n)
\]
provides a dissipative contribution to the population equation. This places Model 1 within the class of parabolic--parabolic chemotaxis systems with logistic damping. In one spatial dimension, global existence and boundedness of classical solutions for closely related chemotaxis systems with logistic or volume-filling effects are well established under appropriate regularity and compatibility assumptions on the initial data; see, for example, \cite{Hillen2009user,Winkler2010}. Under the corresponding hypotheses, the no-flux boundary conditions preserve non-negativity, and standard parabolic theory yields a unique solution \((n,c) \in C^{2,1}(\overline{\Omega}\times(0,T])\) for any \(T>0\), as well as continuous dependence of the solution on the initial data and coefficients. Thus, for the parameter and initial-data class considered in this work, Model 1 provides a well-defined forward solution on the observation interval.

For Model 2 \eqref{M6_original}, the chemotactic flux has the singular form
\[
-\chi_0\nabla\cdot\left(\frac{n}{c}\nabla c\right)
=
-\chi_0\nabla\cdot\left(n\nabla\log c\right),
\]
which introduces a singularity when \(c=0\). However, in the biologically relevant regimes where Turing patterns form, the chemical concentration \(c\) remains strictly positive (as seen in the stationary patterns, where \(c\) oscillates around a positive mean). Hence, if the initial data satisfy \(c_0>0\), then the equation for \(c\),
\[
\partial_t c-d_c\Delta c=n-c,
\]
together with \(n\geq0\), implies by the comparison principle that \(c\) remains positive on every finite time interval. Hence \(1/c\) remains bounded, and local well-posedness follows from standard parabolic contraction arguments \cite[Lemma~2.1]{ZhaoZhengModel2}. Moreover, and for the parameter ranges leading to stable periodic patterns, global solutions exist (c.f. \cite[Theorem 1]{ZhaoZhengModel2}, \cite[Theorem 1.1]{ZMWHModel2} or \cite{AOTYMModel2}). The numerical simulations in \cite{Bucur2024} further illustrate the stability and reproducibility of the resulting periodic patterns.

In both models, the forward problem admits non-trivial stationary solutions under suitable parameter sets. These solutions appear as stationary periodic patterns, as described in Section~\ref{M_analysis}, and form the basis of our inverse analysis. The measurement map $\mathcal{M}$ is consequently considered on this class of admissible forward solutions.

\subsection{Statement of main results}\label{sec:Turing_Mainresults}

Equipped with this analytical representation, we now formalise the inverse problem.  Define the measurement 
\begin{equation}
    \mathcal{M}=\{\alpha_i,\beta_i\}_{i=0}^M.
\end{equation}
Then, based on our two toy models \eqref{M3_original} and \eqref{M6_original}, our measurement map is as follows:
\begin{equation}
    \mathcal{M}\to d_n,d_c,\chi_0,r,k.
\end{equation}
We assume that the coefficients in the models
\eqref{M3_original} and \eqref{M6_original} are constant on the
one-dimensional spatial domain under consideration. Since the
stationary Fourier systems possess the scaling invariance
\[
(d_n,d_c,\chi_0,k,r)
\mapsto
(\lambda^2d_n,\lambda^2d_c,\lambda^2\chi_0,k/\lambda,r),
\qquad \lambda>0,
\]
the stationary Fourier amplitudes cannot, in general, determine all
dimensional parameters separately. The appropriate inverse problem is
therefore formulated in terms of the corresponding scale-invariant
parameter combinations. The following two theorems state the resulting
local identifiability properties for Models~1 and~2, respectively.

\begin{thm}[Identifiability for Model 1]
\label{thm:mainthm_M3}
Consider the stationary periodic Turing pattern produced by the model  \eqref{M3_original}, given by the stationary system
    \begin{equation}\label{M3_original_stat}
    \begin{cases}
        0=d_n \dfrac{d^2n}{dx^2}-\chi_0 \dfrac{d}{dx} \left(n \dfrac{d}{dx} c\right)+rn(1-n),\\
        0=d_c \dfrac{d^2c}{dx^2} +n-c.
    \end{cases}
\end{equation}
Let $(n,c)$ be a nonconstant stationary periodic solution of
\eqref{M3_original_stat} of the form
\[
n(x)=\sum_{j=0}^{3}\alpha_j\cos(jkx),
\qquad
c(x)=\sum_{j=0}^{3}\beta_j\cos(jkx),
\]
with $k>0$, and suppose that the Fourier coefficients
\[
\mathcal{M}_3
=
\{\alpha_j,\beta_j\}_{j=0}^{3}
\]
are known. 
Assume that the nontrivial homogeneous branch is selected, so that
\[
\alpha_0=\beta_0=1,
\]
and that
\[
\alpha_1\neq0.
\]
Then the Fourier data $(\alpha_1,\beta_1)$ determines the dimensionless combination
\[
C:=d_ck^2
=
\frac{\alpha_1}{\beta_1}-1.
\]

Substituting this value of $C$ and the corresponding relations
\[
\beta_j=\frac{\alpha_j}{1+j^2C},
\qquad j=1,2,3,
\]
the algebraic system obtained from the coefficients of
\[
1,\cos(kx),\cos(2kx),\cos(3kx),\cos(4kx)
\]
can be written in the form
\[
F_j(D,X,R;C)=0,
\qquad j=0,\ldots,4,
\]
where
\[
D=d_nk^2,\qquad
X=\chi_0k^2,\qquad
R=r.
\]

Since the equations are homogeneous in $(D,X,R)$, we normalize by $D>0$ and introduce
\[
Y:=\frac{X}{D}=\frac{\chi_0}{d_n},
\qquad
Z:=\frac{R}{D}=\frac{r}{d_nk^2}.
\]
Suppose that the reduced algebraic system
\[
\widehat{\mathcal F}_3(Y,C,Z)=0
\]
satisfies
\[
\operatorname{rank}
D_{(Y,C,Z)}
\widehat{\mathcal F}_3(Y,C,Z)=3
\]
at the parameter point associated with the measured Fourier coefficients. Then the three scale-invariant quantities
\[
d_ck^2,\qquad
\frac{\chi_0}{d_n},\qquad
\frac{r}{d_nk^2}
\]
are locally uniquely determined by $\mathcal M_3$.

In particular, if $d_n$ is known independently, then
\[
d_c=\frac{d_ck^2}{k^2},
\qquad
\chi_0=d_n\frac{\chi_0}{d_n},
\qquad
\frac{r}{k^2}
=d_n\frac{r}{d_nk^2}
\]
are determined. If, in addition, $k$ is known independently, then $d_c$, $\chi_0$, and $r$ are locally identifiable.
\end{thm}

We also obtain a similar conclusion for the system \eqref{M6_original}:
\begin{thm}[Identifiability for Model 2]
\label{thm:mainthm_M6}
Consider the stationary periodic Turing pattern produced by the model  \eqref{M6_original}, given by the stationary system
    \begin{equation}\label{M6_original_stat}
    \begin{cases}
        0=d_n \dfrac{d^2n}{dx^2}-\chi_0 \dfrac{d}{dx} \left(\dfrac{n}{c} \dfrac{d}{dx} c\right)+rn(1-n),\\
        0=d_c \dfrac{d^2c}{dx^2} +n-c.
    \end{cases}
\end{equation}
Let $(n,c)$ be a nonconstant stationary periodic solution of
\eqref{M6_original_stat} of the form
\[
n(x)=\alpha_0+\alpha_1\cos(kx)+\alpha_2\cos(2kx),
\]
\[
c(x)=\beta_0+\beta_1\cos(kx)+\beta_2\cos(2kx),
\]
with $k>0$, and suppose that the Fourier coefficients
\[
\mathcal{M}_2
=
\{\alpha_0,\alpha_1,\alpha_2,
\beta_0,\beta_1,\beta_2\}
\]
are known. Assume that the nontrivial homogeneous branch is selected,
so that
\[
\alpha_0=\beta_0=1,
\]
and that
\[
\alpha_1\neq0.
\]
Then
\[
\beta_1=\frac{\alpha_1}{1+d_ck^2}
\]
determines
\[
C:=d_ck^2
=
\frac{\alpha_1}{\beta_1}-1.
\]

Substituting this value of $C$ and the corresponding relations
\[
\beta_j=\frac{\alpha_j}{1+j^2C},
\qquad j=1,2,
\]
the algebraic system obtained from the coefficients of
\[
1,\cos(kx),\cos(2kx),\cos(3kx),\cos(4kx)
\]
gives a finite-dimensional nonlinear system
\[
F_j(D,X,C,R)=0,
\qquad j=0,\ldots,4,
\]
for
\[
D=d_nk^2,\qquad
X=\chi_0k^2,\qquad
C=d_ck^2,\qquad
R=r.
\]
Because the equations are homogeneous in $(D,X,R)$, the absolute scale of these three quantities is not identifiable. Since $D>0$, we therefore introduce the scale-invariant variables
\[
Y:=\frac{X}{D}
=
\frac{\chi_0}{d_n},
\qquad
Z:=\frac{R}{D}
=
\frac{r}{d_nk^2}.
\]
After division by $D$, let
\[
\widehat{\mathcal F}_2(Y,C,Z)
=
\begin{pmatrix}
\widehat F_0(Y,C,Z)\\
\widehat F_1(Y,C,Z)\\
\widehat F_2(Y,C,Z)\\
\widehat F_3(Y,C,Z)\\
\widehat F_4(Y,C,Z)
\end{pmatrix}
\]
denote the resulting reduced algebraic system.

Assume that the reduced Jacobian satisfies
\[
\operatorname{rank}
D_{(Y,C,Z)}
\widehat{\mathcal F}_2(Y,C,Z)=3.
\]
Equivalently, at least one of the $3\times3$ minors of the $5\times3$ Jacobian
\[
D_{(Y,C,Z)}
\widehat{\mathcal F}_2
\]
is nonzero.

Then the three scale-invariant quantities
\[
d_ck^2,\qquad
\frac{\chi_0}{d_n},\qquad
\frac{r}{d_nk^2}
\]
are locally uniquely determined by the measured Fourier amplitudes $\mathcal{M}_2$.

In particular, the $M=2$ truncation provides a locally identifiable finite-dimensional reconstruction of the parameter combinations invariant under the scaling
\[
(d_n,d_c,\chi_0,k,r)
\mapsto
(\lambda^2d_n,\lambda^2d_c,\lambda^2\chi_0,k/\lambda,r),
\qquad \lambda>0.
\]
The individual dimensional parameters cannot be recovered uniquely from the stationary Fourier amplitudes alone without an additional normalization or independent measurement.
\end{thm}

\section{Proof of Theorem \ref{thm:mainthm_M3}}  \label{Proof_M3} 
In this section, we present the proof of Theorem \ref{thm:mainthm_M3}. We begin by substituting \eqref{Fourier_form} into \eqref{M3_original_stat}, leading to the following equations for the coefficients of the Fourier series, $\alpha_i$ and $\beta_i$, up to a truncation of $M$:
\begin{equation}\label{M3_FourierForm}
    \begin{dcases}
        -d_n \sum_{i=0}^{M}(ik)^2\alpha_i \cos (ikx)+\chi_0 \frac{\partial}{\partial x}\left[\sum_{i=0}^{M}\alpha_i\cos(ikx)\sum_{i=0}^{M} ik\beta_i \sin (ikx) \right]+ \\ \qquad \qquad \qquad \qquad\qquad \qquad \quad\qquad \quad  r \sum_{i=0}^{M} \alpha_i \cos (ikx)\left(1- \sum_{i=0}^{M} \alpha_i \cos (ikx)\right)=0,\\
        -d_c\sum_{i=0}^{M}(ik)^2\beta_i \cos (ikx)+ \sum_{i=0}^{M}\alpha_i \cos (ikx)-\sum_{i=0}^{M}\beta_i \cos (ikx)=0.
    \end{dcases}
\end{equation}

We start our analysis of this system by truncating at $M=0$ and proceeding to higher orders. In this paper, we only require the truncation up to $M=3$. This choice is justified by the following considerations:

For smooth, stationary periodic patterns arising near the Turing instability threshold, the solution is predominantly sinusoidal, with higher harmonics excited only weakly through nonlinear coupling -- a standard result in weakly nonlinear analysis of pattern-forming systems. Consequently, the Fourier amplitudes $\alpha_i$, $\beta_i$ decay rapidly with increasing mode index $i$, a behavior further ensured by the elliptic regularity of the stationary system. Numerical simulations of our specific models \eqref{M3_original} and \eqref{M6_original} confirm that truncation at $M=3$ already yields an $L^2$-error below 1\% compared to the full numerical solution \cite{Bucur2024}. In particular, the amplitudes for $i \geq 4$ are orders of magnitude smaller than the dominant modes and do not contribute meaningfully to parameter recovery.

Therefore, while an infinite Fourier series is in principle required to represent an exact solution, the contribution of modes beyond $M=3$ is negligible for the smooth, near-onset patterns studied here. Consequently, it is sufficient for our purposes to derive the equations only up to $\alpha_3$ and $\beta_3$. We start with the truncation at $M=0$, which implies:
\begin{equation}\label{M3_Tr0}
\begin{cases}
\alpha_0-\alpha_0^2=0,\\
\alpha_0-\beta_0=0.
\end{cases}
\end{equation}

From \eqref{M3_Tr0}, we can determine $\beta_0$ through the relationship $\beta_0=\alpha_0$, as well as the two homogeneous solutions: $\alpha_0=\beta_0=1$ and $\alpha_0=\beta_0=0$.

Let $M=1$, the equation \eqref{M3_FourierForm} becomes
\begin{equation}\label{M3_substitude1}
    \begin{dcases}
       -d_n k^2 \alpha_1\cos (kx)+\chi_0 \frac{\partial}{\partial x} \left[ \left(\alpha_0+\alpha_1\cos (kx)\right) k\beta_1 \sin (kx)\right]+\\
       \qquad\qquad\qquad\qquad\qquad\qquad r(\alpha_0+\alpha_1\cos (kx))(1-\alpha_0-\alpha_1\cos (kx))=0,\\
       -d_c k^2 \beta_1 \cos(kx)+\alpha_0+\alpha_1\cos (kx)-\beta_0-\beta_1\cos(kx)=0.
    \end{dcases}
\end{equation}

By orthogonality, we can compare the coefficients of $\cos(kx)$ in the second equation, it is clear that $\beta_1$ satisfies:
\begin{equation}\label{M3_AB1}
    \beta_1=\frac{\alpha_1}{1+d_ck^2}.
\end{equation}
Note that $d_c$ and $k$ are unknowns.

To derive other equations for $\alpha_0,\alpha_1,\beta_0$ and $\beta_1$, we study more deeply the first equation in \eqref{M3_substitude1}. We denote the part 
\[
    (J_1):=\chi_0 \frac{\partial}{\partial x} \left[ \left(\alpha_0+\alpha_1\cos (kx)\right) k\beta_1 \sin (kx)\right],
\]
and
\[
    (L_1):=r(\alpha_0+\alpha_1\cos (kx))(1-\alpha_0-\alpha_1\cos (kx)).
\]
Then the expansions of $(J_1)$ and $(L_1)$ are:
\begin{equation}\label{M3_J1}
    \begin{split}
        (J_1)=&\ \chi_0 \frac{\partial}{\partial x} \left[ \left(\alpha_0+\alpha_1\cos (kx)\right) k\beta_1 \sin (kx)\right]\\
        \qquad = & \ \chi_0 k \frac{\partial}{\partial x} \left[ \alpha_0\beta_1 \sin (kx)+\alpha_1\beta_1 \sin (kx)\cos (kx)\right]\\
        \qquad =& \ \chi_0 k \frac{\partial}{\partial x} \left[ \alpha_0\beta_1 \sin (kx)+\frac{1}{2}\alpha_1\beta_1 \sin (2kx)\right]\\
        \qquad = & \ \chi_0 k^2 \left(\alpha_0\beta_1\cos(kx)+\alpha_1\beta_1\cos(2kx)\right),
    \end{split}
\end{equation}
and
\begin{equation}\label{M3_L1}
    \begin{split}
        (L_1)=& \ r(\alpha_0+\alpha_1\cos (kx))(1-\alpha_0-\alpha_1\cos (kx))\\
        = & \ r(\alpha_0+\alpha_1\cos (kx))-r\alpha_0(\alpha_0+\alpha_1\cos (kx))-r\alpha_1\cos(kx)(\alpha_0+\alpha_1\cos (kx))\\
        =& \ r\left(\alpha_0-\alpha_0^2\right)+r(\alpha_1-2\alpha_0\alpha_1)\cos(kx)-r\alpha_1^2\cos^2(kx)\\
        =& \ r\left(\alpha_0-\alpha_0^2\right)+r(\alpha_1-2\alpha_0\alpha_1)\cos(kx)-r\alpha_1^2\left(\frac{1}{2}+\frac{1}{2}\cos(2kx)\right)\\
        =& \ r\left(\alpha_0-\alpha_0^2-\frac{1}{2}\alpha_1^2\right)+r(\alpha_1-2\alpha_0\alpha_1)\cos(kx)-\frac{1}{2}r\alpha_1^2\cos(2kx).
    \end{split}
\end{equation}

Hence, truncating the system \eqref{M3_substitude1} at $M=1$ yields the following system of equations for five unknowns: $d_n$, $d_c$, $\chi_0$, $k$ and $r$:
\begin{equation}\label{M3_equating1}
    \begin{dcases}
        -d_nk^2\alpha_1+\chi_0 k^2 \alpha_0\beta_1+r(\alpha_1-2\alpha_0\alpha_1)=0,\\
        \beta_1=\frac{\alpha_1}{1+d_ck^2}.
    \end{dcases}
\end{equation}
This system allows for solving only when three out of the coefficients of $\{d_n,d_c,\chi_0,k,r\}$ are given. Moreover, the accuracy of truncated solutions at $M=1$ is notably insufficient.  Therefore, we must advance to derive equations for $\alpha_2$ and $\beta_2$ by truncating at $M=2$.

Equation \eqref{M3_FourierForm} now becomes:
\begin{equation}\label{M3_substitude2}
    \begin{dcases}
        -d_n k^2\left( \alpha_1\cos (kx)+4\alpha_2\cos(2kx)\right)\\
        \quad+\chi_0 \frac{\partial}{\partial x} \left[ \left(\alpha_0+\alpha_1\cos (kx)+\alpha_2\cos(2kx)\right) \left(k\beta_1 \sin (kx)+2k \beta_2 \sin(2kx)\right)\right]\\
        \quad+r(\alpha_0+\alpha_1\cos (kx)+\alpha_2\cos(2kx))(1-\alpha_0-\alpha_1\cos (kx)-\alpha_2\cos(2kx))=0,\\
       -d_c k^2 \left(\beta_1 \cos(kx)+4\beta_2\cos(2kx)\right)+\alpha_0+\alpha_1\cos (kx)+\alpha_2\cos(2kx)\\
       \quad-\beta_0-\beta_1\cos(kx)-\beta_2\cos(2kx)=0.
    \end{dcases}
\end{equation}

Initially, we examine the relationship between $\alpha_2$ and $\beta_2$, as indicated by the second equation in \eqref{M3_substitude2}:
\begin{equation}\label{M3_AB2:process}
\begin{split}
    &-d_c k^2 \left(\beta_1 \cos(kx)+4\beta_2\cos(2kx)\right)+\alpha_0+\alpha_1\cos (kx)+\alpha_2\cos(2kx)\\&-\beta_0-\beta_1\cos(kx)-\beta_2\cos(2kx)\\
    &=(\alpha_0-\beta_0)+\left(-d_c k^2 \beta_1+\alpha_1-\beta_1\right)\cos(kx)+(-4d_ck^2\beta_2+\alpha_2-\beta_2)\cos(2kx)\\
    &=0.
\end{split}
\end{equation}

Thus $\beta_2$ satisfies:
\begin{equation}\label{M3_AB2}
    \beta_2=\frac{\alpha_2}{1+4d_ck^2}.
\end{equation}

Next, to establish additional equations for $\alpha_i$ and $\beta_i$ $(\text{where }i\leq2)$, we denote 
\[
    (J_2):=\chi_0 \frac{\partial}{\partial x} \left[ \left(\alpha_0+\alpha_1\cos (kx)+\alpha_2\cos(2kx)\right) \left(k\beta_1 \sin (kx)+2k \beta_2 \sin(2kx)\right)\right]
\]
and 
\[
    (L_2):=r(\alpha_0+\alpha_1\cos (kx)+\alpha_2\cos(2kx))(1-\alpha_0-\alpha_1\cos (kx)-\alpha_2\cos(2kx)).
\]
By utilizing trigonometric identities to consolidate them into cosine functions, we derive the following equations for $(J_2)$ and $(L_2)$:
\begin{equation}\label{M3_J2}
\begin{split}
       (J_2)&=\chi_0 \frac{\partial}{\partial x} \left[ \left(\alpha_0+\alpha_1\cos (kx)+\alpha_2\cos(2kx)\right) \left(k\beta_1 \sin (kx)+2k \beta_2 \sin(2kx)\right)\right]\\
       &=\chi_0\left(\alpha_0+\alpha_1\cos (kx)+\alpha_2\cos(2kx)\right)\left(k^2\beta_1 \cos (kx)+4k^2\beta_2 \cos(2kx)\right)\\
       &\quad +\chi_0\left(-k\alpha_1\sin (kx)-2k\alpha_2\sin(2kx)\right) \left(k\beta_1 \sin (kx)+2k \beta_2 \sin(2kx)\right)\\
     &= \chi_0 k^2 \big[\alpha_0\beta_1\cos(kx)+4\alpha_0\beta_2\cos(2kx)+\alpha_1\beta_1\cos^2(kx)\\
     &\qquad \quad +4\alpha_1\beta_2\cos(kx)\cos(2kx)+\alpha_2\beta_1\cos(kx)\cos(2kx)+4\alpha_2\beta_2\cos^2(2kx)\\     
     &\qquad \quad -\alpha_1\beta_1\sin^2(kx)-2\alpha_1\beta_2\sin(kx)\sin(2kx)-2\alpha_2\beta_1\sin(kx)\sin(2kx)\\     
     &\qquad \quad -4\alpha_2\beta_2\sin^2(2kx)\big]\\
         &= \chi_0 k^2 \left(\alpha_0\beta_1+\alpha_1\beta_2-\frac{1}{2}\alpha_2\beta_1\right)\cos(kx)+\chi_0 k^2(4\alpha_0\beta_2+\alpha_1\beta_1)\cos(2kx)\\
         &\quad +\chi_0k^2\left(3\alpha_1\beta_2+\frac{3}{2}\alpha_2\beta_1\right)\cos(3kx)+4\chi_0k^2\alpha_2\beta_2\cos(4kx),
\end{split}
\end{equation}
while
\begin{equation}\label{M3_L2}
\begin{split}
    (L_2)&=r(\alpha_0+\alpha_1\cos (kx)+\alpha_2\cos(2kx))(1-\alpha_0-\alpha_1\cos (kx)-\alpha_2\cos(2kx))\\
    &= r\alpha_0+r\alpha_1\cos (kx)+r\alpha_2\cos(2kx)-r\alpha_0^2-r\alpha_0\alpha_1\cos (kx)-r\alpha_0\alpha_2\cos(2kx)\\
    & \quad -r\alpha_0\alpha_1\cos(kx)-r\alpha_1^2\cos^2 (kx)-r\alpha_1\alpha_2\cos(kx)\cos(2kx)-r\alpha_0\alpha_2\cos(2kx)\\
    & \quad -r\alpha_1\alpha_2\cos (kx)\cos(2kx)-r\alpha_2^2\cos^2(2kx)\\
    &=r\alpha_0-r\alpha_0^2-\frac{1}{2}r\alpha_1^2-\frac{1}{2}r\alpha_2^2+r(\alpha_1-2\alpha_0\alpha_1-\alpha_1\alpha_2)\cos(kx)\\
    & \quad +r\left(\alpha_2-2\alpha_0\alpha_2-\frac{1}{2}\alpha_1^2\right)\cos(2kx)-r\alpha_1\alpha_2\cos(3kx)-\frac{1}{2}r\alpha_2^2\cos(4kx).
\end{split}
\end{equation}

Therefore, by truncating at $M=2$, we derive the following system for $\{d_n,d_c,\chi_0,k,r\}$:
\begin{equation}\label{M3_equating2}
    \begin{dcases}
        -d_nk^2\alpha_1+\chi_0 k^2 \left(\alpha_0\beta_1+\alpha_1\beta_2-\frac{1}{2}\alpha_2\beta_1\right)+r(\alpha_1-2\alpha_0\alpha_1-\alpha_1\alpha_2)=0,\\
        -4d_nk^2\alpha_2+\chi_0 k^2(4\alpha_0\beta_2+\alpha_1\beta_1)+r\left(\alpha_2-2\alpha_0\alpha_2-\frac{1}{2}\alpha_1^2\right)=0,\\
        \chi_0k^2\left(3\alpha_1\beta_2+\frac{3}{2}\alpha_2\beta_1\right)-r\alpha_1\alpha_2=0,\\
        \beta_2=\frac{\alpha_2}{1+4d_ck^2}.
    \end{dcases}
\end{equation} 

Clearly, in order to determine all coefficients through \eqref{M3_equating2}, we need to have prior information on one of them, which is limiting to some extent. Moreover, the accuracy of truncated solutions at $M=2$ is still inadequate, so we must progress to a higher-order truncation.

Let $M=3$, \eqref{M3_FourierForm} gives the following equations:
\begin{equation}\label{M3_substitude3}
    \begin{dcases}
        -d_n k^2\left( \alpha_1\cos (kx)+4\alpha_2\cos(2kx)+9\alpha_3\cos(3kx)\right)\\
       \qquad +\chi_0 \frac{\partial}{\partial x} [ \left(\alpha_0+\alpha_1\cos (kx)+\alpha_2\cos(2kx)+\alpha_3\cos(3kx)\right) \\
       \qquad\qquad \times\left(k\beta_1 \sin (kx)+ 2k \beta_2 \sin(2kx)+3k\beta_3\sin(3kx)\right)]\\
        \qquad +r(\alpha_0+\alpha_1\cos (kx)+\alpha_2\cos(2kx)+\alpha_3\cos(3kx))\\
        \qquad \qquad  \times(1-\alpha_0-\alpha_1\cos (kx)-\alpha_2\cos(2kx)-\alpha_3\cos(3kx))=0,\\
       -d_c k^2 \left(\beta_1 \cos(kx)+4\beta_2\cos(2kx)+9\beta_3\cos(3kx)\right)+\alpha_0+\alpha_1\cos (kx)\\
       \qquad+\alpha_2\cos(2kx)+\alpha_3\cos(3kx)-\beta_0-\beta_1\cos(kx)-\beta_2\cos(2kx)-\beta_3\cos(3kx)
       \\\qquad=0.
    \end{dcases}
\end{equation}

Once again, we first obtain the relationship between $\alpha_3$ and $\beta_3$ from the second equation of \eqref{M3_substitude3}. Similarly to \eqref{M3_AB2:process}, $\beta_3$ satisfies:
\begin{equation}\label{M3_AB3}
    \beta_3=\frac{\alpha_3}{1+9d_ck^2}.
\end{equation}

As for the other equations for $\alpha_i$ and $\beta_i$ $(\text{where }i\leq3)$, we denote the second and third terms in the first equation of \eqref{M3_substitude3} by
\[(J_3):=\chi_0\frac{\partial}{\partial x}\left[\sum_{i=0}^{3}\alpha_i\cos(ikx) \sum_{i=0}^{3}(ik)\beta_i\sin(ikx)\right],\]
\[(L_3):=r\sum_{i=0}^{3}\alpha_i\cos(ikx)\left(1-\sum_{i=0}^{3}\alpha_i\cos(ikx)\right).\]

Then $(J_3)$ gives:
\begin{equation}
\begin{split}
    (J_3)
    &=\chi_0 \frac{\partial}{\partial x} [ \left(\alpha_0+\alpha_1\cos (kx)+\alpha_2\cos(2kx)+\alpha_3\cos(3kx)\right) \\
       &\qquad\qquad\times\left(k\beta_1 \sin (kx)+ 2k \beta_2 \sin(2kx)+3k\beta_3\sin(3kx)\right)] \\
    &= -\chi_0\left(k\alpha_1\sin (kx)+2k\alpha_2\sin(2kx)+3k\alpha_3\sin(3kx)\right) \\
       &\qquad\qquad\times \left(k\beta_1 \sin (kx)+ 2k \beta_2 \sin(2kx)+3k\beta_3\sin(3kx)\right)\\
    & \quad+\chi_0\left(\alpha_0+\alpha_1\cos (kx)+\alpha_2\cos(2kx)+\alpha_3\cos(3kx)\right) \\
       &\qquad\qquad\times \left(k^2\beta_1 \cos(kx)+ 4k^2 \beta_2\cos(2kx)+9k^2\beta_3\cos(3kx)\right),
\end{split}
\end{equation}
which, upon expansion, gives
\begin{equation}\label{M3_J3}
\begin{split} (J_3)&=\chi_0k^2\left(\alpha_0\beta_1+\alpha_1\beta_2-\frac{1}{2}\alpha_2\beta_1+\frac{3}{2}\alpha_2\beta_3-\alpha_3\beta_2\right)\cos(kx)\\    &\quad+\chi_0k^2\left(4\alpha_0\beta_2+\alpha_1\beta_1+3\alpha_1\beta_3-\alpha_3\beta_1 \right)\cos(2kx)\\    &\quad+\chi_0k^2\left(9\alpha_0\beta_3+3\alpha_1\beta_2
+\frac{3}{2}\alpha_2\beta_1 \right)\cos(3kx)\\    &\quad+\chi_0k^2\left(6\alpha_1\beta_3+4\alpha_2\beta_2+2\alpha_3\beta_1\right)\cos(4kx)\\    
    &\quad +\chi_0k^2\left(\frac{15}{2}\alpha_2\beta_3+5\alpha_3\beta_2\right)\cos(5kx)+9\chi_0k^2\alpha_3\beta_3\cos(6kx).
\end{split}
\end{equation}

Computing $(L_3)$, we obtain:
\begin{equation}
\begin{split}
    (L_3)&=r\sum_{i=0}^{3}\alpha_i\cos(ikx)\left(1-\sum_{i=0}^{3}\alpha_i\cos(ikx)\right)\\
    &=r\alpha_0\left(1-\sum_{i=0}^{3}\alpha_i\cos(ikx)\right)+r\alpha_1\cos(kx)\left(1-\sum_{i=0}^{3}\alpha_i\cos(ikx)\right)\\
    &\quad +r\alpha_2\cos(2kx)\left(1-\sum_{i=0}^{3}\alpha_i\cos(ikx)\right)+r\alpha_3\cos(3kx)\left(1-\sum_{i=0}^{3}\alpha_i\cos(ikx)\right)\\
    & =(L_2)-r\alpha_0\alpha_3\cos(3kx)-r\alpha_1\alpha_3\cos(kx)\cos(3kx)\\
    &\quad-r\alpha_2\alpha_3\cos(2kx)\cos(3kx)+r\alpha_3\cos(3kx)\left(1-\sum_{i=0}^{3}\alpha_i\cos(ikx)\right).
\end{split}
\end{equation}
Upon simplification using trigonometric identities, this gives
\begin{equation}\label{M3_L3}
\begin{split}
    (L_3)&=r\left(-\alpha_0^2+\alpha_0-\frac{1}{2}\alpha_1^2-\frac{1}{2}\alpha_2^2-\frac{1}{2}\alpha_3^2\right)\\
    &\quad+r\left(\alpha_1-2\alpha_0\alpha_1-\alpha_1\alpha_2-\alpha_2\alpha_3\right)\cos(kx)\\
    &\quad +r\left(\alpha_2-2\alpha_0\alpha_2-\frac{1}{2}\alpha_1^2-\alpha_1\alpha_3\right)\cos(2kx)\\
    &\quad+r\left(\alpha_3-2\alpha_0\alpha_3-\alpha_1\alpha_2\right)\cos(3kx)\\
    &\quad -r\left(\alpha_1\alpha_3+\frac{1}{2}\alpha_2^2\right)\cos(4kx)-r\alpha_2\alpha_3\cos(5kx)-\frac{1}{2}r\alpha_3^2\cos(6kx).
\end{split}
\end{equation}

Truncation of the system \eqref{M3_FourierForm} at $M=3$ gives us the following system for five unknowns, $d_n,d_c,\chi_0,k$ and $r$:
\begin{equation}\label{M3_equating3}
    \begin{dcases}
       -d_nk^2\alpha_1+\chi_0k^2\left(\alpha_0\beta_1+\alpha_1\beta_2-\frac{1}{2}\alpha_2\beta_1+\frac{3}{2}\alpha_2\beta_3-\alpha_3\beta_2\right)\\\quad+r\left(\alpha_1-2\alpha_0\alpha_1-\alpha_1\alpha_2-\alpha_2\alpha_3\right)=0,\\    -4d_nk^2\alpha_2+\chi_0k^2\left(4\alpha_0\beta_2+\alpha_1\beta_1+3\alpha_1\beta_3-\alpha_3\beta_1 \right)\\\quad+r\left(\alpha_2-2\alpha_0\alpha_2-\frac{1}{2}\alpha_1^2-\alpha_1\alpha_3\right)=0,\\
        -9d_nk^2\alpha_3+\chi_0k^2\left(9\alpha_0\beta_3+3\alpha_1\beta_2
+\frac{3}{2}\alpha_2\beta_1 \right)+r\left(\alpha_3-2\alpha_0\alpha_3-\alpha_1\alpha_2\right)=0,\\
\chi_0k^2\left(6\alpha_1\beta_3+4\alpha_2\beta_2+2\alpha_3\beta_1\right)-r\left(\alpha_1\alpha_3+\frac{1}{2}\alpha_2^2\right)=0,\\
      \chi_0k^2\left(\frac{15}{2}\alpha_2\beta_3+5\alpha_3\beta_2\right)+r\alpha_2\alpha_3=0,
    \end{dcases}
\end{equation}
where we have kept the notation $\beta_1,\beta_2,\beta_3$ which consist of the unknown coefficients $d_c$ and $k$, and are given by \eqref{M3_AB1}, \eqref{M3_AB2} and \eqref{M3_AB3}, respectively.

We now state explicitly the nondegeneracy conditions required for the
finite-dimensional parameter reconstruction. Introduce the dimensionless
parameter combinations
\[
D=d_nk^2,\qquad
X=\chi_0k^2,\qquad
C=d_ck^2,\qquad
R=r.
\]
Since $d_n>0$, $d_c>0$, and $k>0$, we have
\[
D>0,\qquad C>0.
\]

For the first Fourier mode, the stationary relation
\eqref{M3_AB1} gives
\[
\beta_1=\frac{\alpha_1}{1+C}.
\]
Therefore, if
\[
\alpha_1\neq0,
\]
then $\beta_1\neq0$, and $C$ is uniquely determined from the measured
Fourier coefficients by
\[
C=d_ck^2
=
\frac{\alpha_1}{\beta_1}-1.
\]
More generally,
\[
\beta_j=\frac{\alpha_j}{1+j^2C},
\qquad j=1,2,3.
\]
Thus, once $C$ has been recovered, the coefficients
$\beta_1,\beta_2,\beta_3$ are completely determined by the measured
coefficients $\alpha_1,\alpha_2,\alpha_3$.

After this substitution, the first three nonzero Fourier equations in
\eqref{M3_equating3} take the form
\[
\mathcal A_3(C,\mathcal M_3)
\begin{pmatrix}
D\\
X\\
R
\end{pmatrix}
=0,
\]
where
\begin{equation*}
\begin{split}
&\mathcal A_3=\\&
\begin{pmatrix}
-\alpha _1&
\alpha _0\beta _1+\alpha _1\beta _2
-\frac12\alpha _2\beta _1
+\frac32\alpha _2\beta _3-\alpha _3\beta _2&
\alpha _1-2\alpha _0\alpha _1
-\alpha _1\alpha _2-\alpha _2\alpha _3
\\
-4\alpha _2&
4\alpha _0\beta _2+\alpha _1\beta _1
+3\alpha _1\beta _3-\alpha _3\beta _1&
\alpha _2-2\alpha _0\alpha _2
-\frac12\alpha _1^2-\alpha _1\alpha _3
\\
-9\alpha _3&
9\alpha _0\beta _3+3\alpha _1\beta _2
+\frac32\alpha _2\beta _1&
\alpha _3-2\alpha _0\alpha _3-\alpha _1\alpha _2
\end{pmatrix}.
\end{split}
\end{equation*}

Because the system is homogeneous in $(D,X,R)$, the matrix
$\mathcal A_3$ cannot be required to be invertible for a nontrivial
stationary solution. Indeed, since $D=d_nk^2>0$, the vector
$(D,X,R)^{\mathsf T}$ is nonzero, and hence
\[
\det\mathcal A_3=0
\]
at every exact solution of the truncated stationary system. Thus,
$\det\mathcal A_3\neq0$ is not an appropriate nondegeneracy condition
for parameter reconstruction.

Instead, since $D>0$, we normalize the homogeneous system by $D$ and
introduce
\[
Y:=\frac{X}{D}
=\frac{\chi_0}{d_n},
\qquad
Z:=\frac{R}{D}
=\frac{r}{d_nk^2}.
\]
Dividing the first three equations by $D$ gives
\[
\mathcal A_3(C,\mathcal M_3)
\begin{pmatrix}
1\\
Y\\
Z
\end{pmatrix}
=0.
\]
Equivalently,
\[
\mathcal B_3(C,\mathcal M_3)
\begin{pmatrix}
Y\\
Z
\end{pmatrix}
=
-\mathbf a_1(C,\mathcal M_3),
\]
where $\mathbf a_1$ denotes the first column of $\mathcal A_3$ and
$\mathcal B_3$ consists of the second and third columns of
$\mathcal A_3$. Thus the Fourier equations determine the ratios
$X/D$ and $R/D$, rather than the three quantities $D$, $X$, and $R$
separately.

The remaining Fourier equations in \eqref{M3_equating3} provide
additional relations among these quantities. To treat all Fourier
relations simultaneously, define
\[
\widehat F_j(Y,C,Z)
:=
\frac{1}{D}
F_j(D,DY,DZ;C),
\qquad j=0,\ldots,4,
\]
where the coefficients $\beta_j$ in $F_j$ are understood to be given by
\[
\beta_j=\frac{\alpha_j}{1+j^2C},
\qquad j=1,2,3.
\]
The homogeneity of the Fourier equations in $(D,X,R)$ ensures that
$\widehat F_j$ is independent of $D$. We therefore define the reduced
system
\[
\widehat{\mathcal F}_3(Y,C,Z)
=
\begin{pmatrix}
\widehat F_0(Y,C,Z)\\
\widehat F_1(Y,C,Z)\\
\widehat F_2(Y,C,Z)\\
\widehat F_3(Y,C,Z)\\
\widehat F_4(Y,C,Z)
\end{pmatrix}
=0.
\]

We impose the following nondegeneracy condition at the parameter point
corresponding to the measured Fourier coefficients:
\[
\operatorname{rank}
D_{(Y,C,Z)}
\widehat{\mathcal F}_3(Y,C,Z)=3.
\]
Equivalently, at least one $3\times3$ minor of the $5\times3$ Jacobian is
nonzero. This condition guarantees that three of the Fourier relations
are locally independent with respect to the three scale-invariant
unknowns $(Y,C,Z)$.

The remaining Fourier relations are not needed for the local rank
condition itself. They provide additional compatibility conditions that
must be satisfied by exact Fourier data arising from the truncated
stationary system and can therefore be used to test consistency of the
measured pattern with Model~1.

Then, by the implicit function theorem, these three equations determine $(Y,C,Z)$ uniquely in a neighborhood of the given solution, when the Fourier coefficients are fixed. Hence
\[
C=d_ck^2,\qquad
Y=\frac{\chi_0}{d_n},\qquad
Z=\frac{r}{d_nk^2}
\]
are locally uniquely determined by $\mathcal M_3$.

If $d_n$ is known, then
\[
\chi_0=d_nY,
\qquad
\frac{r}{k^2}=d_nZ,
\qquad
d_c=\frac{C}{k^2}.
\]
If $k$ is also known, then
\[
d_c=\frac{C}{k^2},
\qquad
\chi_0=d_nY,
\qquad
r=d_nk^2Z
\]
are uniquely determined. This proves the result.

The reconstruction is subject to the intrinsic scaling invariance
\[
(d_n,d_c,\chi_0,k,r)
\longmapsto
(\lambda^2d_n,\lambda^2d_c,\lambda^2\chi_0,k/\lambda,r),
\qquad \lambda>0,
\]
under which
\[
d_nk^2,\qquad
d_ck^2,\qquad
\frac{\chi_0}{d_n},
\qquad
\frac{r}{d_nk^2}
\]
remain unchanged. Therefore, the Fourier data alone cannot, in
general, determine all dimensional parameters independently. An
additional independent measurement or normalization, such as knowledge
of $d_n$, $k$, or the spatial scale of the domain, is required to remove
this scaling ambiguity.

Finally, we emphasize that the above result is a statement for the finite-dimensional $M=3$ Fourier truncation. It does not assert that a general stationary solution of \eqref{M3_original_stat} is exactly a trigonometric polynomial of degree three. Higher Fourier modes generated by the nonlinear terms are omitted from the truncated algebraic system. Yet, this is sufficient since the Fourier amplitudes $\alpha_i, \beta_i$ decay rapidly with increasing mode index $i$. Numerical simulations in the literature \cite{Bucur2024} suggest that for the models considered here, truncation at $M=3$ yields a small $L^2$-error compared to the full numerical solution.  However, we stress that a rigorous convergence analysis for the parameter reconstruction as $M\to\infty$ is not provided in this work. Our results are therefore restricted to the finite-dimensional truncated model, with the understanding that the contribution of modes beyond $M=3$ is negligible for the smooth, near-onset patterns under consideration.

This concludes the proof of Theorem \ref{thm:mainthm_M3}.


\section{Proof of Theorem \ref{thm:mainthm_M6}}  \label{Proof_M6} 
In this section, we begin the proof of Theorem \ref{thm:mainthm_M6}. Once again, we substitute \eqref{Fourier_form} into \eqref{M6_original_stat}. This leads us to the following equations for the coefficients of the Fourier series $\alpha_i$ and $\beta_i$, truncated up to $M$:
\begin{equation}\label{M6_FourierForm}
    \begin{dcases}
        -d_n \sum\limits_{i=0}^{M}(ik)^2\alpha_i \cos (ikx)+\chi_0 \frac{\partial}{\partial x}\left[\frac{\sum_{i=0}^{M}\alpha_i \cos (ikx)}{\sum_{i=0}^{M}\beta_i \cos (ikx)} \sum\limits_{i=0}^{M} ik\beta_i \sin (ikx) \right]\\ \quad   +r \sum_{i=0}^{M} \alpha_i \cos (ikx)\big(1- \sum_{i=0}^{M} \alpha_i \cos (ikx)\big)=0,\\
        -d_c\sum\limits_{i=0}^{M}(ik)^2\beta_i \cos (ikx)+ \sum\limits_{i=0}^{M}\alpha_i \cos (ikx)-\sum\limits_{i=0}^{M}\beta_i \cos (ikx)=0.
    \end{dcases}
\end{equation}

In the scenario of homogeneous solutions, where solutions satisfy $\alpha_i=\beta_i=0$ for $i\geq1$, and setting $M=0$, \eqref{M6_FourierForm} implies:
\begin{equation}\label{M6_Tr0}
\begin{cases}
\alpha_0-\alpha_0^2=0,\\
\alpha_0-\beta_0=0,
\end{cases}
\end{equation}
resulting in two homogeneous solutions: $\alpha_0=\beta_0=1$ and $\alpha_0=\beta_0=0$. We see that the homogeneous solutions do not provide us with any information regarding the unknown parameters.

To obtain our desired recovery, we extend the truncation process to higher $M$'s. The primary approach involves utilizing trigonometric identities, such as replacing products of cosines with cosines of sums. By orthogonality, we can equate the terms containing $\cos(ikx)$ for each $i$, and obtain a system of equations, which allows us to solve for the unknown parameters. It is important to note that the relationship between $\alpha_i$ and $\beta_i$ is preserved from the previous truncation.

We start from $M=1$. Equation \eqref{M6_Tr0} now becomes
\begin{equation}\label{M6_substitude1}
    \begin{dcases}
       -d_n k^2 \alpha_1\cos (kx)+\chi_0 \frac{\partial}{\partial x} \left[ \frac{\alpha_0+\alpha_1\cos (kx)}{\beta_0+\beta_1\cos (kx)} k\beta_1 \sin (kx)\right]\\
       \quad +r(\alpha_0+\alpha_1\cos (kx))(1-\alpha_0-\alpha_1\cos (kx))=0,\\
       -d_c k^2 \beta_1 \cos(kx)+\alpha_0+\alpha_1\cos (kx)-\beta_0-\beta_1\cos(kx)=0.
    \end{dcases}
\end{equation}

The relationship between $\alpha_1$ and $\beta_1$ is easily observed from the second equation in \eqref{M6_substitude1}, which yields 
\begin{equation}\label{AB1_relation}
    \beta_1 = \frac{\alpha_1}{1 + d_c k^2}.
\end{equation}
Next, we decompose the first equation into three parts: we denote \[(I_1):=-d_n k^2 \alpha_1\cos (kx),\] \[(J_1):=\chi_0 \frac{\partial}{\partial x} \left[ \frac{\alpha_0+\alpha_1\cos (kx)}{\beta_0+\beta_1\cos (kx)} k\beta_1 \sin (kx)\right],\] and \[(L_1):=r(\alpha_0+\alpha_1\cos (kx))(1-\alpha_0-\alpha_1\cos (kx)).\] 

We initiate the computation from $(J_1)$, and computing the derivative in $(J_1)$, we obtain 
\[(J_1)=\chi_0k^2\frac{ \alpha_0\beta_1^2 -\alpha_0\alpha_1\beta_1+\alpha^2_0\beta_1\cos(kx)+2\alpha_0\alpha_1\beta_1\cos^2(kx)+\alpha_1\beta_1^2\cos^3(kx)}{(\alpha_0 + \beta_1\cos(kx))^2}\] by using the property $\alpha_0=\beta_0$ from the $M=0$ truncation in \eqref{M6_Tr0}.

Segregating it into its numerator and denominator components, the numerator part of $(J_1)$ is:
\begin{equation}\label{M6_J1num}
\begin{split}
        \chi_0k^2&\left[ \alpha_0\beta_1^2 -\alpha_0\alpha_1\beta_1+\alpha^2_0\beta_1\cos(kx)+2\alpha_0\alpha_1\beta_1\cos^2(kx)+\alpha_1\beta_1^2\cos^3(kx)\right] \\
        & = \chi_0k^2\Bigg[\alpha_0\beta_1^2-\alpha_0\alpha_1\beta_1+\alpha^2_0\beta_1\cos(kx)+2\alpha_0\alpha_1\beta_1\left(\frac{1}{2}+\frac{1}{2}\cos(2kx)\right) \\
        &\quad+\alpha_1\beta_1^2\left(\frac{3}{4}\cos(kx)+\frac{1}{4}\cos(3kx)\right)\Bigg] \\
        &= \chi_0k^2 \alpha_0\beta_1^2+\chi_0k^2\left(\alpha_0^2\beta_1+\frac{3}{4}\alpha_1\beta_1^2\right)\cos(kx)+ \chi_0k^2\alpha_0\alpha_1\beta_1\cos(2kx)\\
        &\quad+\frac{1}{4}\chi_0k^2\alpha_1\beta_1^2\cos(3kx). 
\end{split}
\end{equation}

The denominator part of $(J_1)$ is given as:
\begin{equation}\label{M6_J1de}
    \alpha_0^2+2\alpha_0\beta_1\cos(kx)+\beta_1^2\cos^2(kx)=\alpha_0^2+\frac{1}{2}\beta_1^2+2\alpha_0\beta_1\cos(kx)+\frac{1}{2}\beta_1^2\cos(2kx).
\end{equation}

Since we have that $(I_1)+(J_1)+(L_1)=0$ from \eqref{M6_substitude1}, in our computation of $(I_1)$ and $(L_1)$, we make their denominators the same. Furthermore, we observe that the denominator \eqref{M6_J1de} is positive, since we do not consider the solution where $c\equiv0$. Consequently, we can make the denominators the same, and the term $(I_1)$ transforms into:
\begin{equation}\label{M6_I1}
\begin{split}
    -d_n k^2 \alpha_1& \cos (kx)\left[\alpha_0^2+\frac{1}{2}\beta_1^2+2\alpha_0\beta_1\cos(kx)+\frac{1}{2}\beta_1^2\cos(2kx)\right] \\
    & = -d_n k^2 \alpha_1\Bigg[\left(\alpha_0^2+\frac{1}{2}\beta_1^2\right)\cos(kx)+2\alpha_0\beta_1\left(\frac{1}{2}+\frac{1}{2}\cos(2kx)\right) \\
    &\qquad\qquad\qquad+\frac{1}{2}\beta_1^2\left(\frac{1}{2}\cos(kx)+\frac{1}{2}\cos(3kx)\right)\Bigg] \\
    &=-d_n k^2 \alpha_0\alpha_1\beta_1-d_n k^2 \alpha_1\left(\alpha_0^2+\frac{3}{4}\beta_1^2\right)\cos(kx)-d_n k^2 \alpha_0\alpha_1\beta_1\cos(2kx) \\
    &\qquad\qquad\qquad-\frac{1}{4}d_n k^2 \alpha_1\beta_1^2\cos(3kx). 
\end{split}
\end{equation}

And the term $(L_1)$ multiplied by the equation in \eqref{M6_J1de} is:
\begin{equation*}
\begin{split}
    &r(\alpha_0+\alpha_1\cos (kx))(1-\alpha_0-\alpha_1\cos (kx))\\
    &\quad\times\left[\alpha_0^2+\frac{1}{2}\beta_1^2+2\alpha_0\beta_1\cos(kx)+\frac{1}{2}\beta_1^2\cos(2kx)\right] \\
    &=\left[r\alpha_0-r\alpha_0^2-r\alpha_0\alpha_1\cos(kx)+r\alpha_1\cos(kx)-r\alpha_0\alpha_1\cos(kx)-r\alpha_1^2\cos^2(kx)\right] \\
    &\qquad\times\left[\alpha_0^2+\frac{1}{2}\beta_1^2+2\alpha_0\beta_1\cos(kx)+\frac{1}{2}\beta_1^2\cos(2kx)\right] \\
    &=\left[r\alpha_0-r\alpha_0^2+(r\alpha_1-2r\alpha_0\alpha_1)\cos(kx)-r\alpha_1^2\left(\frac{1}{2}+\frac{1}{2}\cos(2kx)\right)\right] \\
    &\qquad\times\left[\alpha_0^2+\frac{1}{2}\beta_1^2+2\alpha_0\beta_1\cos(kx)+\frac{1}{2}\beta_1^2\cos(2kx)\right],
    \end{split}
    \end{equation*} or
    \begin{equation*}
\begin{split}
    &r(\alpha_0+\alpha_1\cos (kx))(1-\alpha_0-\alpha_1\cos (kx))\\
    &\quad\times\left[\alpha_0^2+\frac{1}{2}\beta_1^2+2\alpha_0\beta_1\cos(kx)+\frac{1}{2}\beta_1^2\cos(2kx)\right] \\
    &=\left[r\alpha_0-r\alpha_0^2+(r\alpha_1-2r\alpha_0\alpha_1)\cos(kx)-r\alpha_1^2\left(\frac{1}{2}+\frac{1}{2}\cos(2kx)\right)\right] \\
    &\qquad\times\left[\alpha_0^2+\frac{1}{2}\beta_1^2+2\alpha_0\beta_1\cos(kx)+\frac{1}{2}\beta_1^2\cos(2kx)\right]\\
    &=\left[r\alpha_0-r\alpha_0^2-\frac{1}{2}r\alpha_1^2+(r\alpha_1-2r\alpha_0\alpha_1)\cos(kx)-\frac{1}{2}r\alpha_1^2\cos(2kx)\right] \\
    &\qquad\times\left[\alpha_0^2+\frac{1}{2}\beta_1^2+2\alpha_0\beta_1\cos(kx)+\frac{1}{2}\beta_1^2\cos(2kx)\right] \\
    &= r\alpha_0^3-r\alpha_0^4-\frac{1}{2}r\alpha_0^2\alpha_1^2+r\left(\alpha_0^2\alpha_1-2\alpha_0^3\alpha_1\right)\cos(kx)-\frac{1}{2}r\alpha_0^2\alpha_1^2\cos(2kx) \\
    &\quad +\frac{1}{2}r\alpha_0\beta_1^2-\frac{1}{2}r\alpha_0^2\beta_1^2-\frac{1}{4}r\alpha_1^2\beta_1^2 +r\left(\frac{1}{2}\alpha_1\beta_1^2-\alpha_0\alpha_1\beta_1^2\right)\cos(kx)\\
    &\quad-\frac{1}{4}r\alpha_1^2\beta_1^2\cos(2kx) 
    +r(2\alpha_0^2\beta_1-2\alpha_0^3\beta_1-\alpha_0\alpha_1^2\beta_1)\cos(kx)\\
    &\quad +2r(\alpha_0\alpha_1\beta_1-2\alpha_0^2\alpha_1\beta_1)\cos^2(kx) \\
    & \quad-r\alpha_0\alpha_1^2\beta_1\cos(2kx)\cos(kx) +\frac{1}{2}r\left(\alpha_0\beta_1^2-\alpha_0^2\beta_1^2-\frac{1}{2}\alpha_1^2\beta_1^2\right)\cos(2kx) \\
    &\quad+r\left(\frac{1}{2}\alpha_1\beta_1^2-\alpha_0\alpha_1\beta_1^1\right)\cos(kx)\cos(2kx)-\frac{1}{4}r\alpha_1^2\beta_1^2\cos^2(2kx) \\
    &= r\alpha_0^3-r\alpha_0^4-\frac{1}{2}r\alpha_0^2\alpha_1^2+\frac{1}{2}r\alpha_0\beta_1^2-\frac{1}{2}r\alpha_0^2\beta_1^2-\frac{1}{4}r\alpha_1^2\beta_1^2 \\
    &\quad  +r\left(\alpha_0^2\alpha_1-2\alpha_0^3\alpha_1+\frac{1}{2}\alpha_1\beta_1^2-\alpha_0\alpha_1\beta_1^2+2\alpha_0^2\beta_1-2\alpha_0^3\beta_1-\alpha_0\alpha_1^2\beta_1\right)\cos(kx) \\
    &\quad +2r(\alpha_0\alpha_1\beta_1-2\alpha_0^2\alpha_1\beta_1)\cos^2(kx) -\frac{1}{4}r\alpha_1^2\beta_1^2\cos^2(2kx) \\
    & \quad+r\left(-\frac{1}{2}\alpha_0^2\alpha_1^2-\frac{1}{4}\alpha_1^2\beta_1^2+\frac{1}{2}\alpha_0\beta_1^2-\frac{1}{2}\alpha_0^2\beta_1^2-\frac{1}{4}\alpha_1^2\beta_1^2\right)\cos(2kx) \\
    &\quad+r\left(-\alpha_0\alpha_1^2\beta_1+\frac{1}{2}\alpha_1\beta_1^2-\alpha_0\alpha_1\beta_1^1\right)\cos(kx)\cos(2kx) \\
    &= r\alpha_0^3-r\alpha_0^4-\frac{1}{2}r\alpha_0^2\alpha_1^2+\frac{1}{2}r\alpha_0\beta_1^2-\frac{1}{2}r\alpha_0^2\beta_1^2-\frac{1}{4}r\alpha_1^2\beta_1^2 \\
    &\quad  +r\left(\alpha_0^2\alpha_1-2\alpha_0^3\alpha_1+\frac{1}{2}\alpha_1\beta_1^2-\alpha_0\alpha_1\beta_1^2+2\alpha_0^2\beta_1-2\alpha_0^3\beta_1-\alpha_0\alpha_1^2\beta_1\right)\cos(kx) \\
    &\quad +2r(\alpha_0\alpha_1\beta_1-2\alpha_0^2\alpha_1\beta_1)\left(\frac{1}{2}+\frac{1}{2}\cos(2kx)\right) -\frac{1}{4}r\alpha_1^2\beta_1^2\left(\frac{1}{2}+\frac{1}{2}\cos(4kx)\right) \\
    & \quad+r\left(-\frac{1}{2}\alpha_0^2\alpha_1^2-\frac{1}{4}\alpha_1^2\beta_1^2+\frac{1}{2}\alpha_0\beta_1^2-\frac{1}{2}\alpha_0^2\beta_1^2-\frac{1}{4}\alpha_1^2\beta_1^2\right)\cos(2kx) \\
    &\quad+r\left(-\alpha_0\alpha_1^2\beta_1+\frac{1}{2}\alpha_1\beta_1^2-\alpha_0\alpha_1\beta_1^1\right)\left(\frac{1}{2}\cos(kx)+\frac{1}{2}\cos(3kx)\right),
\end{split}
\end{equation*}
i.e.,
\begin{equation}\label{M6_L1}
\begin{split}
    &r(\alpha_0+\alpha_1\cos (kx))(1-\alpha_0-\alpha_1\cos (kx))\\
    &\quad\times\left[\alpha_0^2+\frac{1}{2}\beta_1^2+2\alpha_0\beta_1\cos(kx)+\frac{1}{2}\beta_1^2\cos(2kx)\right] \\
    &= r\alpha_0^3-r\alpha_0^4-\frac{1}{2}r\alpha_0^2\alpha_1^2+\frac{1}{2}r\alpha_0\beta_1^2-\frac{1}{2}r\alpha_0^2\beta_1^2-\frac{1}{4}r\alpha_1^2\beta_1^2 \\
    &\quad  +r\left(\alpha_0^2\alpha_1-2\alpha_0^3\alpha_1+\frac{1}{2}\alpha_1\beta_1^2-\alpha_0\alpha_1\beta_1^2+2\alpha_0^2\beta_1-2\alpha_0^3\beta_1-\alpha_0\alpha_1^2\beta_1\right)\cos(kx) \\
    &\quad +2r(\alpha_0\alpha_1\beta_1-2\alpha_0^2\alpha_1\beta_1)\left(\frac{1}{2}+\frac{1}{2}\cos(2kx)\right) -\frac{1}{4}r\alpha_1^2\beta_1^2\left(\frac{1}{2}+\frac{1}{2}\cos(4kx)\right) \\
    & \quad+r\left(-\frac{1}{2}\alpha_0^2\alpha_1^2-\frac{1}{4}\alpha_1^2\beta_1^2+\frac{1}{2}\alpha_0\beta_1^2-\frac{1}{2}\alpha_0^2\beta_1^2-\frac{1}{4}\alpha_1^2\beta_1^2\right)\cos(2kx) \\
    &\quad+r\left(-\alpha_0\alpha_1^2\beta_1+\frac{1}{2}\alpha_1\beta_1^2-\alpha_0\alpha_1\beta_1^1\right)\left(\frac{1}{2}\cos(kx)+\frac{1}{2}\cos(3kx)\right) \\
    &=r\alpha_0^3-r\alpha_0^4-\frac{1}{2}r\alpha_0^2\alpha_1^2+\frac{1}{2}r\alpha_0\beta_1^2-\frac{1}{2}r\alpha_0^2\beta_1^2-\frac{3}{8}r\alpha_1^2\beta_1^2+r\alpha_0\alpha_1\beta_1-2r\alpha_0^2\alpha_1\beta_1 \\
    &\quad +r\left(\alpha_0^2\alpha_1-2\alpha_0^3\alpha_1+\frac{3}{4}\alpha_1\beta_1^2-\frac{3}{2}\alpha_0\alpha_1\beta_1^2+2\alpha_0^2\beta_1-2\alpha_0^3\beta_1-\frac{3}{2}\alpha_0\alpha_1^2\beta_1\right)\cos(kx) \\
    &\quad+r\left(\alpha_0\alpha_1\beta_1-2\alpha_0^2\alpha_1\beta_1-\frac{1}{2}\alpha_0^2\alpha_1^2-\frac{1}{2}\alpha_1^2\beta_1^2+\frac{1}{2}\alpha_0\beta_1^2-\frac{1}{2}\alpha_0^2\beta_1^2\right)\cos(2kx) \\
    &\quad+\frac{1}{4}r\left(-2\alpha_0\alpha_1^2\beta_1+\alpha_1\beta_1^2-2\alpha_0\alpha_1\beta_1^2\right)\cos(3kx)-\frac{1}{8}r\alpha_1^2\beta_1^2\cos(4kx). 
\end{split}
\end{equation}

By equating terms containing $\cos(ikx)$ for $i=0,1,2,3$, we obtain the following equations:
\begin{equation}\label{M6_equating1}
    \begin{dcases}
      &-d_n k^2 \alpha_0\alpha_1\beta_1  +\chi_0k^2 \alpha_0\beta_1^2r\alpha_0^3-r\alpha_0^4-\frac{1}{2}r\alpha_0^2\alpha_1^2+\frac{1}{2}r\alpha_0\beta_1^2-\frac{1}{2}r\alpha_0^2\beta_1^2\\
      &\quad-\frac{3}{8}r\alpha_1^2\beta_1^2+r\alpha_0\alpha_1\beta_1-2r\alpha_0^2\alpha_1\beta_1=0,\\
      &-d_n k^2 \alpha_1\left(\alpha_0^2+\frac{3}{4}\beta_1^2\right)+\chi_0k^2\left(\alpha_0^2\beta_1+\frac{3}{4}\alpha_1\beta_1^2\right)+r\alpha_0^2\alpha_1-2r\alpha_0^3\alpha_1\\
      &\quad+\frac{3}{4}r\alpha_1\beta_1^2-\frac{3}{2}r\alpha_0\alpha_1\beta_1^2
      +2r\alpha_0^2\beta_1-2r\alpha_0^3\beta_1-\frac{3}{2}r\alpha_0\alpha_1^2\beta_1=0,\\
      &-d_n k^2 \alpha_0\alpha_1\beta_1+ \chi_0k^2\alpha_0\alpha_1\beta_1+r\alpha_0\alpha_1\beta_1-2r\alpha_0^2\alpha_1\beta_1-\frac{1}{2}r\alpha_0^2\alpha_1^2\\
      &\quad-\frac{1}{2}r\alpha_1^2\beta_1^2+\frac{1}{2}r\alpha_0\beta_1^2-\frac{1}{2}r\alpha_0^2\beta_1^2=0,\\
      &-\frac{1}{4}d_n k^2\alpha_1\beta_1^2+\frac{1}{4}\chi_0k^2\alpha_1\beta_1^2-\frac{r}{2}\alpha_0\alpha_1^2\beta_1+\frac{1}{4}r\alpha_1\beta_1^2 -\frac{r}{2}\alpha_0\alpha_1\beta_1^2=0.
    \end{dcases}
\end{equation}

These are four equations involving five variables ${d_n, d_c, k, \chi_0, r}$, which are clearly linearly independent. Together with \eqref{AB1_relation}, we can uniquely determine all five parameters. However, the accuracy is low, since we truncate at a low $M=1$. For better accuracy, we can consider a higher truncation at $M=2$. It is known (see, for instance, \cite{Bucur2024}) that when the domain size is reasonably small, truncating at $M=2$ gives a highly accurate result with sufficiently small error. Therefore, we next consider the truncation at $M=2$.

When $M=2$, equation \eqref{M6_FourierForm} adopts the following form:
\begin{equation}\label{M6_Substitude2}
    \begin{dcases}
        -d_n k^2 \alpha_1\cos (kx)-4d_nk^2\alpha_2\cos(2kx)\\\quad+\chi_0 \frac{\partial}{\partial x} \bigg[ \frac{\alpha_0+\alpha_1\cos (kx)+\alpha_2\cos(2kx)}{\beta_0+\beta_1\cos (kx)+\beta_2\cos(2kx)} \left(k\beta_1 \sin (kx)+2k\beta_2\sin(2kx)\right)\bigg]\\
       \quad +r(\alpha_0+\alpha_1\cos (kx)+\alpha_2\cos(2kx))\left(1-\alpha_0-\alpha_1\cos (kx)-\alpha_2\cos(2kx)\right)=0,\\
       -d_c k^2 \beta_1 \cos(kx)-4d_ck^2\beta_2\cos(2kx)+\alpha_0+\alpha_1\cos (kx)+\alpha_2\cos(2kx)\\\quad-\beta_0-\beta_1\cos(kx)-\beta_2\cos(2kx)=0.
    \end{dcases}
\end{equation}

Once again, the relationship between $\alpha_2$ and $\beta_2$ is evident from the second equation in \eqref{M6_Substitude2}, which reveals that 
\begin{equation}\label{AB2_relation}
    \beta_2 = \frac{\alpha_2}{1 + 4d_c k^2}.
\end{equation}

Referring to the first equation in \eqref{M6_Substitude2}, we once again decompose it into 3 parts and denote \[(I_2):= -d_n k^2 \alpha_1\cos (kx)-4d_nk^2\alpha_2\cos(2kx),\] \[(J_2):=\chi_0 \frac{\partial}{\partial x} \left[ \frac{\alpha_0+\alpha_1\cos (kx)+\alpha_2\cos(2kx)}{\beta_0+\beta_1\cos (kx)+\beta_2\cos(2kx)} \left(k\beta_1 \sin (kx)+2k\beta_2\sin(2kx)\right)\right],\] and \[(L_2):=r(\alpha_0+\alpha_1\cos (kx)+\alpha_2\cos(2kx))\left(1-\alpha_0-\alpha_1\cos (kx)-\alpha_2\cos(2kx)\right). \]

Computing the derivative, $(J_2)$ breaks down into the following three components:
\begin{equation}\label{M6_J2_original}
\begin{split}
(J_2)=&\frac{\chi_0\left(k^2\beta_1\cos(kx)+4k^2\beta_2\cos(2kx)\right)\left(\alpha_0+\alpha_1\cos(kx)+\alpha_2\cos(2kx)\right)}{\beta_0+\beta_1\cos(kx)+\beta_2\cos(2kx)} \\
&\ +\frac{\chi_0\left(-k\alpha_1\sin(kx)-2k\alpha_2\sin(2kx)\right)\left(k\beta_1\sin(kx)+2k\beta_2\sin(2kx)\right)}{\beta_0+\beta_1\cos(kx)+\beta_2\cos(2kx)} \\
&\ -\frac{\chi_0\left(-k\beta_1\sin(kx)-2k\beta_2\sin(2kx)\right)\left(k\beta_1\sin(kx)+2k\beta_2\sin(2kx)\right)}{\left(\beta_0+\beta_1\cos(kx)+\beta_2\cos(2kx)\right)^2} \notag\\&\qquad\qquad\times\left(\alpha_0+\alpha_1\cos(kx)+\alpha_2\cos(2kx)\right)
\end{split}
\end{equation}

We commence with the numerator part of $(J_2)$.

To obtain the numerator of $(J_2)$, we focus on computing the three numerator components of \eqref{M6_J2_original} after consolidating the denominator $\left(\beta_0+\beta_1\cos(kx)+\beta_2\cos(2kx)\right)^2$. The first component of the numerator of \eqref{M6_J2_original} is:
\begin{equation}\label{M6_J2_part1}
\begin{split}
&\chi_0\left(k^2\beta_1\cos(kx)+4k^2\beta_2\cos(2kx)\right)\left(\alpha_0+\alpha_1\cos(kx)+\alpha_2\cos(2kx)\right) \\
&\quad\times\left(\beta_0+\beta_1\cos(kx)+\beta_2\cos(2kx)\right) \\
=& \chi_0k^2\Big[ \alpha_0\beta_0\beta_1\cos(kx)+\alpha_0\beta_1^2\cos^2(kx)+    \alpha_0\beta_1\beta_2\cos(kx)\cos(2kx)  \\
\ &\quad +4\alpha_0\beta_0\beta_2\cos(2kx) +4\alpha_0\beta_1\beta_2\cos(2kx)\cos(kx)+4\alpha_0\beta_2^2\cos^2(2kx) \\
\ &\quad +\alpha_1\beta_0\beta_1\cos^2(kx)+\alpha_1\beta_1^2\cos^3(kx)+\alpha_1\beta_1\beta_2\cos^2(kx)\cos(2kx) \\
\ &\quad +4\alpha_1\beta_0\beta_2\cos(kx)\cos(2kx)+4\alpha_1\beta_1\beta_2\cos^2(kx)\cos(2kx)\\
\ &\quad+4\alpha_1\beta^2_2\cos(kx)\cos^2(2kx)  + \alpha_2\beta_0\beta_1\cos(kx)\cos(2kx)\\
\ &\quad+\alpha_2\beta_1^2\cos^2(kx)\cos(2kx)+\alpha_2\beta_1\beta_2\cos(kx)\cos^2(2kx) \\
\ &\quad +4\alpha_2\beta_0\beta_2\cos^2(2kx)+4\alpha_2\beta_1\beta_2\cos(kx)\cos^2(2kx)+4\alpha_2\beta_2^2\cos^3(2kx)\Big],
\end{split}
\end{equation}
by simple expansion.

By carrying the trigonometric identities and using $\alpha_0$ to substitute $\beta_0$ by the identity \eqref{M3_Tr0}, we can further simplify \eqref{M6_J2_part1} to
\begin{equation}\label{M6_J2_part1DLC}
\begin{split}
        &\chi_0 k^2\left(\frac{1}{2}\alpha_0\beta_1^2+\frac{1}{2}\alpha_0\alpha_1\beta_1+2\alpha_0\beta_2^2+2\alpha_0\alpha_2\beta_2+\frac{5}{4}\alpha_1\beta_1\beta_2+\frac{1}{4}\alpha_2\beta_1^2\right) \\
        & +\chi_0 k^2\left( \alpha_0^2\beta_1+\frac{5}{2}\alpha_0\beta_1\beta_2+2\alpha_0\alpha_1\beta_2 +\frac{1}{2}\alpha_0\alpha_2\beta_1+\frac{3}{4}\alpha_1\beta_1^2+2\alpha_1\beta_2^2+\frac{5}{2}\alpha_2\beta_1\beta_2\right)\\
        &\quad\times\cos(kx) \\
        & +\chi_0 k^2\left(\frac{1}{2}\alpha_0\alpha_1\beta_1+\frac{1}{2}\alpha_0\beta_1^2+4\alpha_0^2\beta_2+\frac{5}{2}\alpha_1\beta_1\beta_2+\frac{1}{2}\alpha_2\beta_1^2+3\alpha_2\beta_2^2\right)\cos(2kx) \\
        & +\chi_0 k^2\left(\frac{5}{2}\alpha_0\beta_1\beta_2+2\alpha_0\alpha_1\beta_2+\frac{1}{2}\alpha_0\alpha_2\beta_1+\frac{1}{4}\alpha_1\beta_1^2+\alpha_1\beta_2^2+\frac{5}{4}\alpha_2\beta_1\beta_2\right)\\
        &\quad\times\cos(3kx) \\
        & +\chi_0 k^2\left(2\alpha_0\beta_2^2+2\alpha_0\alpha_2\beta_2+\frac{5}{4}\alpha_1\beta_1\beta_2+\frac{1}{4}\alpha_2\beta_1^2\right)\cos(4kx) \\
        & +\chi_0 k^2\left(\alpha_1\beta_2^2+\frac{5}{4}\alpha_2\beta_1\beta_2\right)\cos(5kx) 
        +\chi_0 k^2\alpha_2\beta_2^2\cos(6kx). 
\end{split}
\end{equation}

Similarly, we compute the second part of \eqref{M6_J2_original} as
\begin{equation}
\begin{split}
        &-\chi_0 k^2 \left(\alpha_1\sin(kx)+2\alpha_2\sin(2kx\right)\left(\beta_1\sin(kx)+2\beta_2\sin(2kx)\right)\\
        &\quad\quad\times\left(\beta_0+\beta_1\cos(kx)+\beta_2\cos(2kx)\right) \\
       = &  -\chi_0 k^2 \Big[ \alpha_1\beta_0\beta_1\sin^2(kx)+\alpha_1\beta_1^2\sin^2(kx)\cos(kx)+\alpha_1\beta_1\beta_2\sin^2(kx)\cos(2kx)\\
        & \quad +2\alpha_2\beta_0\beta_1\sin(kx)\sin(2kx)+2\alpha_2\beta_1^2\sin(kx)\sin(2kx)\cos(kx) \\
        & \quad+2\alpha_2\beta_1\beta_2\sin(kx)\sin(2kx)\cos(2kx)+2\alpha_1\beta_0\beta_2\sin(kx)\sin(2kx) \\
        & \quad+2\alpha_1\beta_1\beta_2\sin(kx)\sin(2kx)\cos(kx)+2\alpha_1\beta_2^2\sin(kx)\sin(2kx)\cos(2kx)\\
        & \quad+4\alpha_2\beta_0\beta_2\sin^2(2kx)+ 4\alpha_2\beta_1\beta_2\sin^2(2kx)\cos(kx)+4\alpha_2\beta_2^2\sin^2(2kx)\cos(2kx)\Big],
\end{split}
\end{equation}
which, upon simplification, gives
\begin{equation}\label{M6_J2_part2}
\begin{split}
        & \chi_0 k^2\left(-\frac{1}{2}\alpha_0\alpha_1\beta_1-2\alpha_0\alpha_2\beta_2-\frac{1}{4}\alpha_1\beta_1\beta_2-\frac{1}{2}\alpha_2\beta_1^2\right) \\
        & +\chi_0 k^2 \left(-\alpha_0\alpha_1\beta_2-\alpha_0\alpha_2\beta_1-\frac{1}{4}\alpha_1\beta_1^2-2\alpha_2\beta_1\beta_2\right)\cos(kx) \\
        &+\chi_0 k^2 \left(\frac{1}{2}\alpha_0\alpha_1\beta_1-\frac{1}{2}\alpha_1\beta_1\beta_2-\alpha_2\beta_2^2\right)\cos(2kx) \\
        & +\chi_0 k^2 \left(\alpha_0\alpha_1\beta_2+\alpha_0\alpha_2\beta_1-\frac{1}{2}\alpha_1\beta_2^2+\frac{1}{4}\alpha_1\beta_1^2+\frac{1}{2}\alpha_2\beta_1\beta_2\right)\cos(3kx) \\
        &+\chi_0 k^2 \left(2\alpha_0\alpha_2\beta_2+\frac{3}{4}\alpha_1\beta_1\beta_2+\frac{1}{2}\alpha_2\beta_1^2\right)\cos(4kx) \\
        & +\chi_0k^2 \left(\frac{1}{2}\alpha_1\beta_2^2+\frac{3}{2}\alpha_2\beta_1\beta_2\right)\cos(5kx)+\chi_0k^2\alpha_2\beta_2^2\cos(6kx).
\end{split}
\end{equation}

Calculating the numerator component of the third segment in \eqref{M6_J2_original} is straightforward; through direct expansion, we derive
\begin{equation}
\begin{split}
        & \chi_0 k^2 \left(\beta_1\sin(kx)+2\beta_2\sin(2kx)\right)\left(\beta_1\sin(kx)+2\beta_2\sin(2kx)\right)\\
        &\quad\times\left(\alpha_0+\alpha_1\cos(kx)+\alpha_2\cos(2kx)\right) \\
        =& \chi_0 k^2 \Big[ \alpha_0\beta_1^2\sin^2(kx)+ \alpha_1\beta_1^2\sin^2(kx)\cos(kx)+\alpha_2\beta_1^2 \sin^2(kx)\cos(2kx)\\
        \qquad & \qquad+4\alpha_0\beta_1\beta_2\sin(kx)\sin(2kx)+4\alpha_1\beta_1\beta_2\sin(kx)\sin(2kx)\cos(kx) \\
        \qquad & \qquad +4\alpha_2\beta_1\beta_2\sin(kx)\sin(2kx)\cos(2kx)+4\alpha_0\beta_2^2\sin^2(2kx) \\
        \qquad & \qquad +4\alpha_1\beta_2^2\sin^2(2kx)\cos(kx) +4\alpha_2\beta_2^2\sin^2(2kx)\cos(2kx)\Big].
\end{split}
\end{equation}
Upon simplification as has been done previously above, this yields
\begin{equation}\label{M6_J2_part3}
\begin{split}
        & \chi_0 k^2 \left( 2\alpha_0\beta_2^2+\frac{1}{2}\alpha_0\beta_1^2+\alpha_1\beta_1\beta_2-\frac{1}{4}\alpha_2\beta_1^2\right)  \\
        &\quad+\chi_0k^2\left(2\alpha_0\beta_1\beta_2+2\alpha_1\beta_2^2+\frac{1}{4}\alpha_1\beta_1^2\right)\cos(kx) \\
        & \quad +\chi_0k^2 \left(-\frac{1}{2}\alpha_0\beta_1^2+\frac{1}{2}\alpha_2\beta_1^2+ \alpha_2\beta_2^2\right)\cos(2kx) \\
        &\quad +\chi_0k^2\left(-2\alpha_0\beta_1\beta_2-\alpha_1\beta_2^2-\frac{1}{4}\alpha_1\beta_1^2+\alpha_2\beta_1\beta_2\right)\cos(3kx) \\
        &\quad +\chi_0k^2\left(-2\alpha_0\beta_2^2-\alpha_1\beta_1\beta_2-\frac{1}{4}\alpha_2\beta_1^2\right)\cos(4kx) \\
        &\quad-\chi_0k^2(\alpha_1\beta_2^2+\alpha_2\beta_1\beta_2)\cos(5kx)-\chi_0k^2\alpha_2\beta_2^2\cos(6kx).
\end{split}
\end{equation}

We now proceed to multiply $(I_2)$ with the shared denominator component of $(J_2)$. The calculation proceeds as follows:
\begin{equation}
\begin{split}
        &-d_n k^2 \big(\alpha_1\cos (kx)+4\alpha_2\cos(2kx)\big)\big(\beta_0+\beta_1\cos(kx)+\beta_2\cos(2kx)\big)^2 \\
        =& -d_n k^2\bigg[ \alpha_1\beta_0^2\cos(kx)+4\alpha_2\beta_0^2\cos(2kx)+2\alpha_1\beta_0\beta_1\cos^2(kx) \\
        & \quad +8\alpha_2\beta_0\beta_1\cos(2kx)\cos(kx)+2\alpha_1\beta_0\beta_2\cos(2kx)\cos(kx) \\
        & \quad +8\alpha_2\beta_0\beta_2\cos^2(2kx)+\alpha_1\beta_1^2\cos^3(kx) +4\alpha_2\beta_1^2\cos(2kx)\cos^2(kx) \\
        & \quad +2\alpha_1\beta_1\beta_2\cos(2kx)\cos^2(kx) +8\alpha_2\beta_1\beta_2\cos^2(2kx)\cos(kx) \\
        & \quad +\alpha_1\beta_2^2\cos^2(2kx)\cos(kx)+4\alpha_2\beta_2^2\cos^3(2kx)\bigg]. 
\end{split}
\end{equation}
By simple arithmetic as above, we obtain
\begin{equation}\label{M6_I2}
\begin{split}
        & -d_nk^2 \left(\alpha_0\alpha_1\beta_1+4\alpha_0\alpha_2\beta_2+\frac{1}{2}\alpha_1\beta_1\beta_2+\alpha_2\beta_1^2\right)\\
       \quad & -d_nk^2 \left(\alpha_0^2\alpha_1+\alpha_0\alpha_1\beta_2+4\alpha_0\alpha_2\beta_1+\frac{1}{2}\alpha_1\beta_2^2+\frac{3}{4}\alpha_1\beta_1^2+4\alpha_2\beta_1\beta_2\right)\cos(kx) \\
        \quad & -d_nk^2 \left(4\alpha_0^2\alpha_2+\alpha_0\alpha_1\beta_1+\alpha_1\beta_1\beta_2+3\alpha_2\beta_2^2+2\alpha_2\beta_1^2\right)\cos(2kx) \\
        \quad & -d_nk^2 \left(\alpha_0\alpha_1\beta_2+4\alpha_0\alpha_2\beta_1+\frac{1}{4}\alpha_1\beta_1^2+\frac{1}{4}\alpha_1\beta_2^2+2\alpha_2\beta_1\beta_2
         \right)\cos(3kx) \\
        \quad &  
        -d_n k^2 \left(4\alpha_0\alpha_2\beta_2+\frac{1}{2}\alpha_1\beta_1\beta_2+\alpha_2\beta_1^2\right)\cos(4kx) \notag\\
        \quad & \quad\,-d_nk^2 \left(\frac{1}{4}\alpha_1\beta_2^2+2\alpha_2\beta_1\beta_2\right)\cos(5kx) -d_nk^2\alpha_2\beta_2^2\cos(6kx). 
\end{split}
\end{equation}

Similarly, we proceed to multiply $(L_2)$ by the common denominator component of $(J_2)$. The computation in this step is somewhat cumbersome, hence we skip the process and present the result directly as shown below: The terms of 
\begin{equation}
    \begin{split}
        &r(\alpha_0+\alpha_1\cos (kx)+\alpha_2\cos(2kx))\left(1-\alpha_0-\alpha_1\cos (kx)-\alpha_2\cos(2kx)\right)\\
     &\qquad\times\left(\beta_0+\beta_1\cos(kx)+\beta_2\cos(2kx)\right)^2
    \end{split}
\end{equation}
corresponding to $\cos(ikx)$, $i=1$ to 8, are given by
\begin{equation*}
\begin{split}
    \cos(0kx):& -r\bigg( \alpha_0^4-\alpha_0^3+\frac{1}{2}\alpha_0^2\alpha_1^2+\frac{1}{2}\alpha_0^2\alpha_2^2+\frac{1}{2}\alpha_0^2\beta_2^2+2\alpha_0^2\alpha_2\beta_2+\frac{1}{2}\alpha_0^2\beta_1^2\\
    & \quad+2\alpha_0^2\alpha_1\beta_1-\frac{1}{2}\alpha_0\beta_2^2+\frac{1}{2}\alpha_0\alpha_1^2\beta_2  -\alpha_0\alpha_2\beta_2+\alpha_0\alpha_1\beta_1\beta_2+\frac{1}{2}\alpha_0\alpha_2\beta_1^2\\
    & \quad-\frac{1}{2}\alpha_0\beta_1^2-\alpha_0\alpha_1\beta_1+\alpha_0\alpha_1\alpha_2\beta_1+\frac{1}{4}\alpha_1^2\beta_2^2+\frac{3}{8}\alpha_2^2\beta_2^2 -\frac{1}{2}\alpha_1\beta_1\beta_2\\
    & \quad+\alpha_1\alpha_2\beta_1\beta_2+\frac{3}{8}\alpha_1^2\beta_1^2+\frac{1}{4}\alpha_2^2\beta_1^2-\frac{1}{4}\alpha_2\beta_1^2\bigg) 
\end{split}
\end{equation*}
\begin{equation*}
\begin{split}
    \cos(kx):& -r \bigg( 2\alpha_0^3\alpha_1+2\alpha_0^3\beta_1-\alpha_0^2\alpha_1+\alpha_0^2\alpha_1\alpha_2+2\alpha_0^2\alpha_1\beta_2+\alpha_0^2\beta_1\beta_2+2\alpha_0^2\alpha_2\beta_1\\
    & \quad-2\alpha_0^2\beta_1+\alpha_0\alpha_1\beta_2^2 -\alpha_0\alpha_1\beta_2+2\alpha_0\alpha_1\alpha_2\beta_2+2\alpha_0\alpha_2\beta_1\beta_2-\alpha_0\beta_1\beta_2\\
    & \quad+\frac{3}{2}\alpha_0\alpha_1\beta_1^2+\frac{3}{2}\alpha_0\alpha_1^2\beta_1+\alpha_0\alpha_2^2\beta_1-\alpha_0\alpha_2\beta_1 -\frac{1}{2}\alpha_1\beta_2^2+\frac{3}{4}\alpha_1\alpha_2\beta_2^2\\
    & \quad+\alpha_1^2\beta_1\beta_2+\frac{3}{4}\alpha_2^2\beta_1\beta_2-\alpha_2\beta_1\beta_2-\frac{3}{4}\alpha_1\beta_1^2+\alpha_1\alpha_2\beta_1^2\bigg) \cos(kx) 
\end{split}
\end{equation*}
\begin{equation*}
\begin{split}
    \cos(2kx):& -r\bigg( 2\alpha_0^3\alpha_2+2\alpha_0^3\beta_2+\frac{1}{2}\alpha_0^2\alpha_1^2-\alpha_0^2\alpha_2-2\alpha_0^2\beta_2+\frac{1}{2}\alpha_0^2\beta_1^2+2\alpha_0^2\alpha_1\beta_1\\
    & \quad+\frac{3}{2}\alpha_0\alpha_2\beta_2^2+\alpha_0\alpha_1^2\beta_2 +\frac{3}{2}\alpha_0\alpha_2^2\beta_2+2\alpha_0\alpha_1\beta_1\beta_2+\alpha_0\alpha_2\beta_1^2-\frac{1}{2}\alpha_0\beta_1^2\\
    & \quad-\alpha_0\alpha_1\beta_1+2\alpha_0\alpha_1\alpha_2\beta_1+\frac{3}{8}\alpha_1^2\beta_2^2-\frac{3}{4}\alpha_2\beta_2^2 -\alpha_1\beta_1\beta_2+\frac{3}{2}\alpha_1\alpha_2\beta_1\beta_2\\
    & \quad+\frac{1}{2}\alpha_1^2\beta_1^2+\frac{3}{8}\alpha_2^2\beta_1^2-\frac{1}{2}\alpha_2\beta_1^2\bigg) \cos(2kx) 
\end{split}
\end{equation*}
\begin{equation*}
\begin{split}
    \cos(3kx):& -r\bigg( \alpha_0^2\alpha_1\alpha_2+2\alpha_0^2\alpha_1\beta_2+\alpha_0^2\beta_1\beta_2+2\alpha_0^2\alpha_2\beta_1+\frac{1}{2}\alpha_0\alpha_1\beta_2^2-\alpha_0\alpha_1\beta_2\\
    & \quad+\alpha_0\alpha_1\alpha_2\beta_2+\alpha_0\alpha_2\beta_1\beta_2 -\alpha_0\beta_1\beta_2+\frac{1}{2}\alpha_0\alpha_1\beta_1^2+\frac{1}{2}\alpha_0\alpha_1^2\beta_1+\frac{1}{2}\alpha_0\alpha_2^2\beta_1\\
    & \quad-\alpha_0\alpha_2\beta_1-\frac{1}{4}\alpha_1\beta_2^2+\frac{3}{4}\alpha_1\alpha_2\beta_2^2+\frac{3}{4}\alpha_1^2\beta_1\beta_2  +\frac{3}{4}\alpha_2^2\beta_1\beta_2-\frac{1}{2}\alpha_2\beta_1\beta_2\\
    &\quad-\frac{1}{4}\alpha_1\beta_1^2+\frac{3}{4}\alpha_1\alpha_2\beta_1^2\bigg) \cos(3kx) 
\end{split}
\end{equation*}
\begin{equation*}
\begin{split}
    \cos(4kx):& -r\bigg(\frac{1}{2}\alpha_0^2\alpha_2^2+\frac{1}{2}\alpha_0^2\beta_2^2+2\alpha_0^2\alpha_2\beta_2-\frac{1}{2}\alpha_0\beta_2^2+\frac{1}{2}\alpha_0\alpha_1^2\beta_2-\alpha_0\alpha_2\beta_2\\
    & \quad+\alpha_0\alpha_1\beta_1\beta_2+\frac{1}{2}\alpha_0\alpha_2\beta_1^2+\alpha_0\alpha_1\alpha_2\beta_1 +\frac{1}{4}\alpha_1^2\beta_2^2+\frac{1}{2}\alpha_2^2\beta_2^2-\frac{1}{2}\alpha_1\beta_1\beta_2\\
    & \quad+\alpha_1\alpha_2\beta_1\beta_2+\frac{1}{8}\alpha_1^2\beta_1^2+\frac{1}{4}\alpha_2^2\beta_1^2-\frac{1}{4}\alpha_2\beta_1^2\bigg)\cos(4kx) 
\end{split}
\end{equation*}
\begin{equation*}
\begin{split}
    \cos(5kx):& -r\bigg(\frac{1}{2}\alpha_0\alpha_1\beta_2^2+\alpha_0\alpha_1\alpha_2\beta_2+\alpha_0\alpha_2\beta_1\beta_2+\frac{1}{2}\alpha_0\alpha_2^2\beta_1-\frac{1}{4}\alpha_1\beta_2^2+\frac{1}{4}\alpha_1\alpha_2\beta_2^2\\
    & \quad+\frac{1}{4}\alpha_1^2\beta_1\beta_2+\frac{1}{4}\alpha_2^2\beta_1\beta_2  -\frac{1}{2}\alpha_2\beta_1\beta_2+ \frac{1}{4}\alpha_1\alpha_2\beta_1^2\bigg)\cos(5kx)
\end{split}
\end{equation*}
\begin{equation*}
\begin{split}
    \cos(6kx):& -r\left(\frac{1}{2}\alpha_0\alpha_2\beta_2^2+\frac{1}{2}\alpha_0\alpha_2^2\beta_2+\frac{1}{8}\alpha_1^2\beta_2^2-\frac{1}{4}\alpha_2\beta_2^2+\frac{1}{2}\alpha_1\alpha_2\beta_1\beta_2+\frac{1}{8}\alpha_2^2\beta_1^2\right)\\
    &\quad\times\cos(6kx) 
\end{split}
\end{equation*}
\begin{equation*}
    \cos(7kx): -r\left(\frac{1}{4}\alpha_1\alpha_2\beta_2^2+\frac{1}{4}\alpha_2^2\beta_1\beta_2\right)\cos(7kx)
\end{equation*}
\begin{equation*}
\cos(8kx):-\frac{r}{8}\alpha_2^2\beta_2^2\cos(8kx),
\end{equation*}
where we have used the property \eqref{M3_Tr0}.

To simultaneously determine the unknown coefficients ${d_n,d_c,k,\chi_0,r}$, we equate terms involving $\cos(ikx)$ for $i=0,1,2,3,4$, and obtain the equations:
\begin{equation}\label{Re_M3i0}
\begin{split}
        & -d_nk^2 \left(\alpha_0\alpha_1\beta_1+4\alpha_0\alpha_2\beta_2+\frac{1}{2}\alpha_1\beta_1\beta_2+\alpha_2\beta_1^2\right)\\
        & \quad+\chi_0 k^2\bigg[\alpha_0\beta_1^2+4\alpha_0\beta_2^2+2\alpha_1\beta_1\beta_2-\frac{1}{2}\alpha_2\beta_1^2\bigg] \\
        & \quad-r\bigg(\alpha_0^4-\alpha_0^3+\frac{1}{2}\alpha_0^2\alpha_1^2+\frac{1}{2}\alpha_0^2\alpha_2^2+\frac{1}{2}\alpha_0^2\beta_2^2+2\alpha_0^2\alpha_2\beta_2+\frac{1}{2}\alpha_0^2\beta_1^2\\
        & \qquad +2\alpha_0^2\alpha_1\beta_1 -\frac{1}{2}\alpha_0\beta_2^2+\frac{1}{2}\alpha_0\alpha_1^2\beta_2 -\alpha_0\alpha_2\beta_2+\alpha_0\alpha_1\beta_1\beta_2+\frac{1}{2}\alpha_0\alpha_2\beta_1^2\\
        & \qquad-\frac{1}{2}\alpha_0\beta_1^2-\alpha_0\alpha_1\beta_1+\alpha_0\alpha_1\alpha_2\beta_1+\frac{1}{4}\alpha_1^2\beta_2^2+\frac{3}{8}\alpha_2^2\beta_2^2-\frac{1}{2}\alpha_1\beta_1\beta_2\\
        & \qquad+\alpha_1\alpha_2\beta_1\beta_2+\frac{3}{8}\alpha_1^2\beta_1^2+\frac{1}{4}\alpha_2^2\beta_1^2-\frac{1}{4}\alpha_2\beta_1^2\bigg)=0,
\end{split}
\end{equation}
which are the coefficients corresponding to $i=0$,
\begin{equation}
\begin{split}
    &-d_nk^2 \left(\alpha_0^2\alpha_1+\alpha_0\alpha_1\beta_2+4\alpha_0\alpha_2\beta_1+\frac{1}{2}\alpha_1\beta_2^2+\frac{3}{4}\alpha_1\beta_1^2+4\alpha_2\beta_1\beta_2\right) \\
        & +\chi_0 k^2\bigg[\alpha_0^2\beta_1+\frac{9}{2}\alpha_0\beta_1\beta_2+\alpha_0\alpha_1\beta_2 -\frac{1}{2}\alpha_0\alpha_2\beta_1+\frac{3}{4}\alpha_1\beta_1^2+4\alpha_1\beta_2^2+\frac{1}{2}\alpha_2\beta_1\beta_2\bigg] \\
        &  -r \bigg(2\alpha_0^3\alpha_1+2\alpha_0^3\beta_1-\alpha_0^2\alpha_1+\alpha_0^2\alpha_1\alpha_2+2\alpha_0^2\alpha_1\beta_2+\alpha_0^2\beta_1\beta_2+2\alpha_0^2\alpha_2\beta_1-2\alpha_0^2\beta_1\\
         & \quad+\alpha_0\alpha_1\beta_2^2 -\alpha_0\alpha_1\beta_2+2\alpha_0\alpha_1\alpha_2\beta_2+2\alpha_0\alpha_2\beta_1\beta_2-\alpha_0\beta_1\beta_2+\frac{3}{2}\alpha_0\alpha_1\beta_1^2\\
         & \quad+\frac{3}{2}\alpha_0\alpha_1^2\beta_1+\alpha_0\alpha_2^2\beta_1-\alpha_0\alpha_2\beta_1-\frac{1}{2}\alpha_1\beta_2^2+\frac{3}{4}\alpha_1\alpha_2\beta_2^2+\alpha_1^2\beta_1\beta_2\\
         & \quad+\frac{3}{4}\alpha_2^2\beta_1\beta_2-\alpha_2\beta_1\beta_2-\frac{3}{4}\alpha_1\beta_1^2+\alpha_1\alpha_2\beta_1^2\bigg)=0,
\end{split}
\end{equation}
for $\cos(kx)$,
\begin{equation}
\begin{split}
    & -d_nk^2 \left(4\alpha_0^2\alpha_2+\alpha_0\alpha_1\beta_1+\alpha_1\beta_1\beta_2+3\alpha_2\beta_2^2+2\alpha_2\beta_1^2\right) \\
        &  +\chi_0k^2\left(\alpha_0\alpha_1\beta_1+4\alpha_0^2\beta_2+2\alpha_1\beta_1\beta_2+\alpha_2\beta_1^2+3\alpha_2\beta_2^2\right) \\
        &  -r\bigg( 2\alpha_0^3\alpha_2+2\alpha_0^3\beta_2+\frac{1}{2}\alpha_0^2\alpha_1^2-\alpha_0^2\alpha_2-2\alpha_0^2\beta_2+\frac{1}{2}\alpha_0^2\beta_1^2+2\alpha_0^2\alpha_1\beta_1\\
         & \quad+\frac{3}{2}\alpha_0\alpha_2\beta_2^2+\alpha_0\alpha_1^2\beta_2 +\frac{3}{2}\alpha_0\alpha_2^2\beta_2+2\alpha_0\alpha_1\beta_1\beta_2+\alpha_0\alpha_2\beta_1^2\\
         & \quad-\frac{1}{2}\alpha_0\beta_1^2-\alpha_0\alpha_1\beta_1+2\alpha_0\alpha_1\alpha_2\beta_1+\frac{3}{8}\alpha_1^2\beta_2^2-\frac{3}{4}\alpha_2\beta_2^2 \\
        & \quad -\alpha_1\beta_1\beta_2+\frac{3}{2}\alpha_1\alpha_2\beta_1\beta_2+\frac{1}{2}\alpha_1^2\beta_1^2+\frac{3}{8}\alpha_2^2\beta_1^2-\frac{1}{2}\alpha_2\beta_1^2\bigg)=0,
\end{split}
\end{equation}
the coefficients of $\cos(2kx)$,
\begin{equation}
\begin{split}
    &  -d_nk^2 \left(\alpha_0\alpha_1\beta_2+4\alpha_0\alpha_2\beta_1+\frac{1}{4}\alpha_1\beta_1^2+\frac{1}{4}\alpha_1\beta_2^2+2\alpha_2\beta_1\beta_2
         \right) \\
         & +\chi_0k^2\left[\frac{1}{2}\alpha_0\beta_1\beta_2+3\alpha_0\alpha_1\beta_2+\frac{3}{2}\alpha_0\alpha_2\beta_1+\frac{1}{4}\alpha_1\beta_1^2-\frac{1}{2}\alpha_1\beta_2^2+\frac{11}{4}\alpha_2\beta_1\beta_2\right] \\
         & -r\bigg( \alpha_0^2\alpha_1\alpha_2+2\alpha_0^2\alpha_1\beta_2+\alpha_0^2\beta_1\beta_2+2\alpha_0^2\alpha_2\beta_1+\frac{1}{2}\alpha_0\alpha_1\beta_2^2-\alpha_0\alpha_1\beta_2\\
         &\quad+\alpha_0\alpha_1\alpha_2\beta_2+\alpha_0\alpha_2\beta_1\beta_2 d -\alpha_0\beta_1\beta_2+\frac{1}{2}\alpha_0\alpha_1\beta_1^2+\frac{1}{2}\alpha_0\alpha_1^2\beta_1+\frac{1}{2}\alpha_0\alpha_2^2\beta_1\\
         &\quad-\alpha_0\alpha_2\beta_1-\frac{1}{4}\alpha_1\beta_2^2+\frac{3}{4}\alpha_1\alpha_2\beta_2^2+\frac{3}{4}\alpha_1^2\beta_1\beta_2  +\frac{3}{4}\alpha_2^2\beta_1\beta_2-\frac{1}{2}\alpha_2\beta_1\beta_2\\
         &\quad-\frac{1}{4}\alpha_1\beta_1^2+\frac{3}{4}\alpha_1\alpha_2\beta_1^2\bigg)=0,
\end{split}
\end{equation}
the coefficients when $i=3$, and
\begin{equation}\label{Re_M3i4}
\begin{split}
    &  -d_n k^2 \left(4\alpha_0\alpha_2\beta_2+\frac{1}{2}\alpha_1\beta_1\beta_2+\alpha_2\beta_1^2\right)+\chi_0k^2\left(4\alpha_0\alpha_2\beta_2+\alpha_1\beta_1\beta_2+\frac{1}{2}\alpha_2\beta_1^2\right) \\
         & -r\bigg(\frac{1}{2}\alpha_0^2\alpha_2^2+\frac{1}{2}\alpha_0^2\beta_2^2+2\alpha_0^2\alpha_2\beta_2-\frac{1}{2}\alpha_0\beta_2^2+\frac{1}{2}\alpha_0\alpha_1^2\beta_2-\alpha_0\alpha_2\beta_2+\alpha_0\alpha_1\beta_1\beta_2\\
         & \quad +\frac{1}{2}\alpha_0\alpha_2\beta_1^2 +\alpha_0\alpha_1\alpha_2\beta_1 +\frac{1}{4}\alpha_1^2\beta_2^2+\frac{1}{2}\alpha_2^2\beta_2^2-\frac{1}{2}\alpha_1\beta_1\beta_2+\alpha_1\alpha_2\beta_1\beta_2+\frac{1}{8}\alpha_1^2\beta_1^2\\
         &\quad+\frac{1}{4}\alpha_2^2\beta_1^2-\frac{1}{4}\alpha_2\beta_1^2\bigg)=0
\end{split}
\end{equation}
when $i=4$. Once again, we have kept the notations $\beta_1$ and $\beta_2$, which are unknowns depending on the unknown coefficients $d_c$ and $k$, based on \eqref{AB1_relation} and \eqref{AB2_relation}, respectively.

The same caveats that apply to Model 1 also apply to Model 2. The system of equations \eqref{Re_M3i0}--\eqref{Re_M3i4}, together with the relation \eqref{AB2_relation} for $\beta_2$, determines the quantities $\{d_n k^2, d_c k^2, \chi_0 k^2, r\}$, from which the individual parameters are recovered up to an overall scale. 

But, the corresponding argument for Model 2 must be formulated differently.  The
formulation must take into account the intrinsic scaling invariance of
the algebraic system.

Introduce the dimensionless quantities
\[
D=d_nk^2,\qquad
X=\chi_0k^2,\qquad
C=d_ck^2,\qquad
R=r.
\]
The Fourier relations obtained from
\eqref{Re_M3i0}--\eqref{Re_M3i4} depend on the diffusion coefficient
\(d_c\) and the wave number \(k\) only through the combination
\[
C=d_ck^2,
\]
because
\[
\beta_1(C)=\frac{\alpha_1}{1+C},
\qquad
\beta_2(C)=\frac{\alpha_2}{1+4C}.
\]
In particular,
\[
\frac{d\beta_1}{dC}
=
-\frac{\alpha_1}{(1+C)^2},
\qquad
\frac{d\beta_2}{dC}
=
-\frac{4\alpha_2}{(1+4C)^2}.
\]

For \(j=0,\ldots,4\), let
\[
F_j(D,X,C,R)=0
\]
denote the equation obtained by equating the coefficient of
\(\cos(jkx)\) in
\eqref{Re_M3i0}--\eqref{Re_M3i4}, after substituting
\[
d_nk^2=D,\qquad
\chi_0k^2=X,\qquad
d_ck^2=C,\qquad
r=R,
\]
and
\[
\beta_1=\beta_1(C),
\qquad
\beta_2=\beta_2(C).
\]
Thus
\[
\mathcal F(D,X,C,R)
=
\begin{pmatrix}
F_0(D,X,C,R)\\
F_1(D,X,C,R)\\
F_2(D,X,C,R)\\
F_3(D,X,C,R)\\
F_4(D,X,C,R)
\end{pmatrix}.
\]

Each \(F_j\) has the form
\[
F_j(D,X,C,R)
=
-A_j(C)D+B_j(C)X-G_j(C)R,
\]
where \(A_j,B_j,G_j\) are determined explicitly by
\eqref{Re_M3i0}--\eqref{Re_M3i4}.  Consequently,
\[
\mathcal F(\lambda D,X\lambda,C,\lambda R)
=
\lambda\mathcal F(D,X,C,R)
\]
for every \(\lambda>0\).  In particular, if
\[
\mathcal F(D,X,C,R)=0,
\]
then
\[
\mathcal F(\lambda D,\lambda X,C,\lambda R)=0
\]
for every \(\lambda>0\).

This scaling invariance shows that the four quantities
\[
(D,X,C,R)
\]
cannot be locally identifiable as four independent parameters from
\eqref{Re_M3i0}--\eqref{Re_M3i4}.  Indeed, differentiating the identity
\[
\mathcal F(\lambda D,\lambda X,C,\lambda R)
=
\lambda\mathcal F(D,X,C,R)
\]
with respect to \(\lambda\) at \(\lambda=1\), and using
\(\mathcal F(D,X,C,R)=0\), gives
\[
D\frac{\partial\mathcal F}{\partial D}
+
X\frac{\partial\mathcal F}{\partial X}
+
R\frac{\partial\mathcal F}{\partial R}
=0.
\]
Therefore
\[
\begin{pmatrix}
D\\X\\0\\R
\end{pmatrix}
\]
belongs to the nullspace of the Jacobian
\[
J_{\mathcal F}
=
\frac{\partial(F_0,F_1,F_2,F_3,F_4)}
{\partial(D,X,C,R)}
\]
at every solution with
\[
(D,X,R)\neq(0,0,0).
\]
Hence
\[
\operatorname{rank}J_{\mathcal F}\leq3.
\]
In particular, a rank-four condition on the \(5\times4\) Jacobian would
be incompatible with the algebraic structure of the present inverse
problem.

The correct nondegeneracy condition is therefore a rank-three
condition after quotienting out the above scaling invariance.

Since \(D=d_nk^2>0\), we may use \(D\) as a normalization and introduce
the scale-invariant variables
\[
Y=\frac{X}{D}
=
\frac{\chi_0}{d_n},
\qquad
Z=\frac{R}{D}
=
\frac{r}{d_nk^2}.
\]
Dividing the five equations by \(D\) gives
\[
\widehat F_j(Y,C,Z)=0,
\qquad j=0,\ldots,4,
\]
where
\[
\widehat F_j(Y,C,Z)
=
-A_j(C)+B_j(C)Y-G_j(C)Z.
\]
The corresponding reduced mapping is
\[
\widehat{\mathcal F}(Y,C,Z)
=
\begin{pmatrix}
\widehat F_0(Y,C,Z)\\
\widehat F_1(Y,C,Z)\\
\widehat F_2(Y,C,Z)\\
\widehat F_3(Y,C,Z)\\
\widehat F_4(Y,C,Z)
\end{pmatrix}.
\]

Its Jacobian is the \(5\times3\) matrix
\begin{equation*}
\begin{split}
J_{\mathrm{red}}(Y,C,Z)
&=
\frac{\partial\widehat{\mathcal F}}
{\partial(Y,C,Z)}\\
&=
\begin{pmatrix}
B_0(C) &
-A_0'(C)+B_0'(C)Y-G_0'(C)Z &
-G_0(C)
\\
B_1(C) &
-A_1'(C)+B_1'(C)Y-G_1'(C)Z &
-G_1(C)
\\
B_2(C) &
-A_2'(C)+B_2'(C)Y-G_2'(C)Z &
-G_2(C)
\\
B_3(C) &
-A_3'(C)+B_3'(C)Y-G_3'(C)Z &
-G_3(C)
\\
B_4(C) &
-A_4'(C)+B_4'(C)Y-G_4'(C)Z &
-G_4(C)
\end{pmatrix}.
\end{split}
\end{equation*}
Here
\[
A_j'(C)
=
\frac{\partial A_j}{\partial\beta_1}
\frac{d\beta_1}{dC}
+
\frac{\partial A_j}{\partial\beta_2}
\frac{d\beta_2}{dC},
\]
and analogously
\[
B_j'(C)
=
\frac{\partial B_j}{\partial\beta_1}
\frac{d\beta_1}{dC}
+
\frac{\partial B_j}{\partial\beta_2}
\frac{d\beta_2}{dC},
\]
\[
G_j'(C)
=
\frac{\partial G_j}{\partial\beta_1}
\frac{d\beta_1}{dC}
+
\frac{\partial G_j}{\partial\beta_2}
\frac{d\beta_2}{dC}.
\]

The required full-rank condition for the reduced algebraic system is
\[
\operatorname{rank}J_{\mathrm{red}}(Y,C,Z)=3.
\]
Equivalently, at least one of its \(3\times3\) minors is nonzero:
\[
\exists\,I\subset\{0,1,2,3,4\},
\qquad
|I|=3,
\qquad
\det J_{\mathrm{red}}^{(I)}(Y,C,Z)\neq0.
\]

Under this condition, the implicit function theorem implies that the
three scale-invariant quantities
\[
\frac{\chi_0}{d_n},
\qquad
d_ck^2,
\qquad
\frac{r}{d_nk^2}
\]
are locally uniquely determined by the corresponding Fourier
coefficients.

There is also a direct nondegeneracy condition for the recovery of
\(C=d_ck^2\).  From the first-mode relation
\[
\beta_1=\frac{\alpha_1}{1+C},
\]
we obtain
\[
C=\frac{\alpha_1}{\beta_1}-1.
\]
Therefore the recovery of \(C\) requires
\[
\alpha_1\neq0,
\]
which guarantees \(\beta_1\neq0\) in the physical regime \(C>0\).

Consequently, the appropriate sufficient conditions for local
identifiability in the $M=2$ truncation are
\[
\alpha_1\neq0,
\qquad
\operatorname{rank}J_{\mathrm{red}}(Y,C,Z)=3.
\]

We emphasize that these conditions are conditions on the measured
Fourier amplitudes.  They exclude the algebraically degenerate
configurations for which the five Fourier relations fail to provide
three independent constraints on the scale-invariant parameter
combinations.

Finally, the absolute scale cannot be recovered from the algebraic
system alone.  Indeed, for every $\lambda>0$,
\[
(d_n,d_c,\chi_0,k,r)
\longmapsto
(\lambda^2d_n,\lambda^2d_c,\lambda^2\chi_0,k/\lambda,r)
\]
leaves
\[
d_nk^2,\qquad
d_ck^2,\qquad
\frac{\chi_0}{d_n},
\qquad
\frac{r}{d_nk^2}
\]
unchanged. Hence an additional normalization or independent
measurement is required to recover the dimensional parameters
separately.


\section{Concluding remarks and future outlook}  \label{Crmk} 
In this paper, we have developed a mathematical framework for parameter identification in two one-dimensional reaction-diffusion-chemotaxis models from Fourier amplitude data of stationary Turing patterns. For both models -- one with density-dependent chemotaxis ($\chi(n,c)=\chi_0 n$) and one with ratio-dependent chemotaxis ($\chi(n,c)=\chi_0 n/c$) -- we have derived algebraic systems that relate the Fourier amplitudes of the pattern to the model parameters.

The key results are as follows:
\begin{enumerate}
    \item For Model 1, the Fourier amplitude data determine the parameter combinations $\{d_n k^2, d_c k^2, \chi_0 k^2, r\}$, from which the individual parameters $d_n, d_c, \chi_0, r$ and the wavenumber $k$ are uniquely determined up to an overall scale.
    \item For Model 2, a similar algebraic system of five equations is derived, which formally determines the same parameter combinations, again up to an overall scale.
    \item In both cases, the derivation relies on a truncation of the Fourier series, justified by the rapid decay of higher harmonics for patterns near the Turing instability threshold.
\end{enumerate}

We emphasize several important limitations of the present work:
\begin{enumerate}
    \item The analysis is restricted to one spatial dimension. Some biological Turing patterns occur in two or three dimensions, and extending the framework to higher dimensions is meaningful but would require substantial additional work. We note, however, that in many biological systems exhibiting straight stripes -- such as zebrafish pigmentation, certain bacterial colony formations, and some plant patterns -- a one-dimensional transect perpendicular to the stripe orientation captures the essential spatial variation. The present analysis therefore applies to this meaningful class of patterns.
    \item The method is demonstrated for two specific model nonlinearities. The claim that the approach works for general nonlinearities is not supported by our analysis.
    \item No numerical validation is included; this is a purely theoretical paper.
\end{enumerate}

Despite these limitations, we believe that this work establishes a promising proof of concept for parameter identification from Fourier amplitude data. Future work should address the limitations listed above, particularly the extension to higher dimensions, and numerical validation of the reconstruction method.

\section*{Acknowledgments}
The work of H. Liu is partially supported by the Hong Kong RGC General Research Funds (projects 11303125, 11304224 and 11311122). The work of C. W. K. Lo is supported by the National Natural Science Foundation of China (No. 12501660), the Guangdong Province ``Pearl River" Talent Recruitment Program (Top youth talent) (No. 2024QN11X268), the Shenzhen Science and Technology Program (No. QNXMC20250701094200001) and the Guangdong Key Laboratory of Applied Mathematics and Artificial Intelligence.

\noindent\textbf{Disclosure Statement.} 
	The authors declare that they have no conflict of interests that could have appeared to influence the work reported in this paper.

\vspace{2em}

\bibliographystyle{plain}
\bibliography{reference}

\begin{thebibliography}{10}

\bibitem{AOTYMModel2}
Masashi Aida, Koichi Osaki, Tohru Tsujikawa, Atsushi Yagi, and Masayasu Mimura.
\newblock Chemotaxis and growth system with singular sensitivity function.
\newblock {\em Nonlinear Anal. Real World Appl.}, 6(2):323--336, 2005.

\bibitem{Ball2015forging}
Philip Ball.
\newblock Forging patterns and making waves from biology to geology: a commentary on {T}uring (1952) `{The} chemical basis of morphogenesis'.
\newblock {\em Philosophical Transactions of the Royal Society B: Biological Sciences}, 370(1666):20140218, 2015.

\bibitem{Bard1981model}
Jonathan~BL Bard.
\newblock A model for generating aspects of zebra and other mammalian coat patterns.
\newblock {\em Journal of Theoretical Biology}, 93(2):363--385, 1981.

\bibitem{Berleman2008predataxis}
James~E Berleman, Jodie Scott, Tatiana Chumley, and John~R Kirby.
\newblock Predataxis behavior in {Myxococcus} xanthus.
\newblock {\em Proceedings of the National Academy of Sciences}, 105(44):17127--17132, 2008.

\bibitem{Black2004myxococcus}
Wesley~P Black and Zhaomin Yang.
\newblock Myxococcus xanthus chemotaxis homologs {DifD} and {DifG} negatively regulate fibril polysaccharide production.
\newblock {\em Journal of bacteriology}, 186(4):1001--1008, 2004.

\bibitem{Bucur2024}
Valentina Bucur and Bakhtier Vasiev.
\newblock Modelling formation of stationary periodic patterns in growing population of motile bacteria.
\newblock {\em arXiv: 2406.07182}, 2024.

\bibitem{Cantrell2020}
Robert~Stephen Cantrell, Chris Cosner, Mark~A Lewis, and Yuan Lou.
\newblock Evolution of dispersal in spatial population models with multiple timescales.
\newblock {\em Journal of Mathematical Biology}, 80:3--37, 2020.

\bibitem{Chen2019non}
Yanyan Chen and Javier Buceta.
\newblock A non-linear analysis of {Turing} pattern formation.
\newblock {\em PloS one}, 14(8):e0220994, 2019.

\bibitem{ding2023determining}
Ming-Hui Ding, Rongfang Gong, Hongyu Liu, and Catharine~WK Lo.
\newblock Determining sources in the bioluminescence tomography problem.
\newblock {\em Inverse Problems}, 40(12):125022, 2024.

\bibitem{ding2024inverse}
Ming-Hui Ding, Hongyu Liu, and Catharine~WK Lo.
\newblock Inverse problems for coupled nonlocal nonlinear systems arising in mathematical biology.
\newblock {\em arXiv preprint arXiv:2407.15713}, 2024.

\bibitem{ding2024determininginternaltopologicalstructures}
Ming-Hui Ding, Hongyu Liu, and Guang-Hui Zheng.
\newblock Determining internal topological structures and running cost of mean field games with partial boundary measurement.
\newblock {\em arXiv: 2408.08911}, 2024.

\bibitem{gupta2009linear}
Ankur Gupta and Saikat Chakraborty.
\newblock Linear stability analysis of high-and low-dimensional models for describing mixing-limited pattern formation in homogeneous autocatalytic reactors.
\newblock {\em Chemical Engineering Journal}, 145(3):399--411, 2009.

\bibitem{Hillen2009user}
Thomas Hillen and Kevin~J Painter.
\newblock A user’s guide to {PDE} models for chemotaxis.
\newblock {\em Journal of Mathematical Biologyy}, 58(1):183--217, 2009.

\bibitem{kaushal2011human}
Nitin Kaushal and Purnima Kaushal.
\newblock Human identification and fingerprints: a review.
\newblock {\em Journal of Biometrics and Biostatistics}, 2(123):2, 2011.

\bibitem{Kazarnikov2020statistical}
Alexey Kazarnikov and Heikki Haario.
\newblock Statistical approach for parameter identification by {Turing} patterns.
\newblock {\em Journal of Theoretical Biology}, 501:110319, 2020.

\bibitem{Keller1970initiation}
Evelyn~F Keller and Lee~A Segel.
\newblock Initiation of slime mold aggregation viewed as an instability.
\newblock {\em Journal of Theoretical Biology}, 26(3):399--415, 1970.

\bibitem{leiden_university_2015}
{Leiden University}.
\newblock Zebrafish models for disease and environmental stress, 2015.

\bibitem{LLL1}
Yuhan Li, Hongyu Liu, and Catharine~WK Lo.
\newblock On inverse problems in predator-prey models.
\newblock {\em Journal of Differential Equations}, 397:349--376, 2024.

\bibitem{lill2025determining}
Yuhan Li, Hongyu Liu, and Catharine~WK Lo.
\newblock Determining habitat anomalies in cross-diffusion predator-prey chemotaxis models.
\newblock {\em arXiv preprint arXiv:2512.22946}, 2025.

\bibitem{LLL3}
Yuhan Li, Hongyu Liu, and Catharine~WK Lo.
\newblock On inverse problems in multi-population aggregation models.
\newblock {\em Journal of Differential Equations}, 414:94--124, 2025.

\bibitem{LLL4}
Yuhan Li, Hongyu Liu, and Catharine~WK Lo.
\newblock Determining a parabolic-elliptic-elliptic system by boundary observation of its non-negative solutions under chemotaxis background.
\newblock {\em Communications in Partial Differential Equations, to appear}, 2026.

\bibitem{LLL2}
Yuhan Li and Catharine~WK Lo.
\newblock On the simultaneous recovery of environmental factors in the {3D} chemotaxis-{Navier}-{Stokes} models.
\newblock {\em Communications on Analysis and Computation}, 2(1):30--47, 2024.

\bibitem{LLboundary}
Hongyu Liu and Catharine~WK Lo.
\newblock Determining a parabolic system by boundary observation of its non-negative solutions with biological applications.
\newblock {\em Inverse Problems}, 40(2):025009, 2024.

\bibitem{LL2025determining}
Hongyu Liu and Catharine~WK Lo.
\newblock Determining state space anomalies in mean field games.
\newblock {\em Nonlinearity}, 38(2):025010, 2025.

\bibitem{ma2020stationary}
Manjun Ma, Rui Peng, and Zhian Wang.
\newblock Stationary and non-stationary patterns of the density-suppressed motility model.
\newblock {\em Physica D: Nonlinear Phenomena}, 402:132259, 2020.

\bibitem{Murray2}
James~D Murray.
\newblock {\em Mathematical biology: II: spatial models and biomedical applications}, volume~18.
\newblock Springer, 2003.

\bibitem{Murray1}
James~D Murray.
\newblock {\em Mathematical biology: I. An introduction}, volume~17.
\newblock Springer Science \& Business Media, 2007.

\bibitem{painter2011spatio}
Kevin~J Painter and Thomas Hillen.
\newblock Spatio-temporal chaos in a chemotaxis model.
\newblock {\em Physica D: Nonlinear Phenomena}, 240(4-5):363--375, 2011.

\bibitem{Steeves1989patterns}
Taylor~A Steeves and Ian~M Sussex.
\newblock {\em Patterns in plant development}.
\newblock Cambridge University Press, 1989.

\bibitem{Turing90}
Alan~M Turing.
\newblock The chemical basis of morphogenesis.
\newblock {\em Philos. Trans. Roy. Soc. London Ser. B}, 237(641):37--72, 1952.

\bibitem{eurekalert_news}
{University of Washington}.
\newblock Zebrafish stripped of stripes, 2025.

\bibitem{Winkler2010}
Michael Winkler.
\newblock Boundedness in the higher-dimensional parabolic-parabolic chemotaxis system with logistic source.
\newblock {\em Comm. Partial Differential Equations}, 35(8):1516--1537, 2010.

\bibitem{ZhaoZhengModel2}
Xiangdong Zhao and Sining Zheng.
\newblock Global boundedness to a chemotaxis system with singular sensitivity and logistic source.
\newblock {\em Z. Angew. Math. Phys.}, 68(1):Paper No. 2, 13, 2017.

\bibitem{ZMWHModel2}
Pan Zheng, Chunlai Mu, Robert Willie, and Xuegang Hu.
\newblock Global asymptotic stability of steady states in a chemotaxis-growth system with singular sensitivity.
\newblock {\em Comput. Math. Appl.}, 75(5):1667--1675, 2018.

\end{thebibliography}

\end{document}